\pdfoutput=1
\documentclass[12pt,a4paper]{article}
\usepackage{float}
\usepackage[english]{babel}
\usepackage{epsf}
\usepackage{latexsym}
\usepackage{epsfig}
\usepackage {amssymb}
\usepackage {amsmath}
\usepackage{setspace}
\usepackage{booktabs}
\usepackage[colorlinks,linkcolor=black, citecolor=black]{hyperref}
\usepackage{enumitem}

\usepackage{rotating}

\usepackage[longnamesfirst]{natbib}
\bibpunct{(}{)}{;}{a}{,}{,}
\renewcommand{\cite}{\citep}
\newcommand{\biblist}{
\bibliographystyle{apalike}
\setlength{\bibsep}{0cm}
}

\newtheorem{theorem}{Theorem}

\newtheorem{lemma}{Lemma}
\newtheorem{corollary}{Corollary}
\newtheorem{proposition}{Proposition}

\def\beqn{\begin{eqnarray*}}
\def\eeqn{\end{eqnarray*}}
\def\beq{\begin{eqnarray}}
\def\eeq{\end{eqnarray}}

\newcommand{\N}{{\mathbb N}}

\newcommand{\R}{{\mathbb R}}
\newcommand{\Z}{{\mathbb Z}}

\newcommand{\eop}{\hfill $\square$\par}
\newcommand{\Prob}{\mathrm{P}}
\newcommand{\diff}{\mathrm{d}}

\newcommand{\LTwo}{\calL^2\left(\Prob^{X} \right)}

\newcommand{\alphaest}{\widehat{\alpha}}
\newcommand{\HS}{\mathrm{HS}}
\newcommand{\NCov}{\Sigma}
\newcommand{\KCov}{S}
\newcommand{\KEst}{S_n}
\newcommand{\bEst}{T_n^\ast y}

\newcommand{\ssize}{.98\linewidth}
\def\T{{ \mathrm{\scriptscriptstyle T} }}
\def \calL {\mathcal L}

\newcommand{\fest}[2]{f_{#1}^{[#2]}}
\newcommand{\pol}[2]{p_{#1}^{[#2]}}
\newcommand{\qpol}[2]{q_{#1}^{[#2]}}
\newcommand{\xr}[3]{x_{#1,#2}^{[#3]}}
\newcommand{\polinner}[3]{\left[#1,#2\right]_{#3}}
\newcommand{\polnorm}[2]{\polinner{#1}{#1}{#2}}
\newcommand{\pdiff}[1]{\left|\left(\pol{#1}{1}\right)'(0)\right|}
\newcommand{\pdifftwo}[2]{\left|\left(\pol{#1}{#2}\right)'(0)\right|}
\newcommand{\ed}{d_\lambda}
\newcommand{\aopt}{{a^\ast}}
\newcommand{\Proj}{P}
\newcommand{\Sg}{\KCov^{1/2}}
\newcommand{\Sgs}{\KCov}
\newcommand{\Ctau}{C}
\newcommand{\E}{\mathrm{E}}
\newcommand{\trace}{\mathrm{tr}}

  \def \calA {\mathcal A}

 \def \calH {\mathcal H}

 \def \calK {\mathcal K}
 \def \calL {\mathcal L}

 \def \calP {\mathcal P}

\makeatother
\makeatletter
\renewcommand{\@fnsymbol}[1]{\@arabic{#1}}
\makeatother

\title{Kernel partial least squares for stationary data}
\author{Marco Singer$^a$,
\hspace{0.3pt} Tatyana Krivobokova$^a$,
\hspace{0.3pt} Axel Munk$^{a,b}$
\medskip\\
$^a$Institute for Mathematical Stochastics, G\"ottingen, Germany
\medskip\\
$^b$Max Planck Institute for Biophysical Chemistry, G\"ottingen, Germany
}

\begin{document}
\baselineskip=15pt
\maketitle

\begin{abstract}
\baselineskip=15pt
We consider the kernel partial least squares algorithm for non-parametric regression with stationary dependent data. Probabilistic convergence rates of the kernel partial least squares estimator to the true regression function are established under a source and an effective dimensionality conditions. It is shown both theoretically and in simulations that long range dependence results in slower convergence rates. A protein dynamics example shows high predictive power of kernel partial least squares.
\vspace{15pt}

\end{abstract}


\pagenumbering{arabic}

\section{Introduction}
\label{sec:Introduction}

Partial least squares (PLS) is a regularized regression technique developed by \citet{Wol84} to deal with collinearities in the regressor matrix. It is an iterative algorithm where the covariance between response and regressor is maximized at each step, see \citet{Hel88} for a detailed description. Regularization in the PLS algorithm is obtained by stopping the iteration process early.

Several studies showed that partial least squares algorithm is competitive with other regression methods such as ridge regression and principal component regression and it needs generally fewer iterations than the latter to achieve comparable estimation and prediction, see, e.g., \citet{FrankFriedman} and \citet{Krae07b}. For an overview of further properties of PLS we refer to \citet{Ros06}.

Reproducing kernel Hilbert spaces (RKHS) have a long history in probability and statistics \citep[see e.g.][]{Berlinet}. 
Here we focus on the supervised kernel based learning approach for the solution of non-parametric regression problems. RKHS methods are both computationally and theoretically attractive, due to the kernel trick \citep{Sch98} and the representer theorem \citep{Wah99} as well as its generalization \citep{Sch01}. Within the reproducing kernel Hilbert space framework one can adapt linear regularized regression techniques like ridge regression and principal component regression to a non-parametric setting, see \citet{Sau98} and \citet{Ros00a}, respectively. We refer to \citet{bSch} for more details on the kernel based learning approach.

Kernel PLS was introduced in \citet{Ros01} who reformulated the algorithm presented in \citet{Lin93}. The relationship to kernel conjugate gradient (KCG) methods was highlighted in \citet{Blan10a}. It can be seen in \citet{Hanke} that conjugate gradient methods are well suited for handling ill-posed problems, as they arise in kernel learning, see, e.g., \citet{Vit06}.

\citet{Ros03} investigated the performance of kernel partial least squares (KPLS) for non-linear discriminant analysis.
\citet{Blan10a} proved the consistency of KPLS when the algorithm is stopped early without giving convergence rates.  

\citet{Cap07} showed that kernel ridge regression (KRR) attains optimal probabilistic rates of convergence for independent and identically distributed data, using a source and a polynomial effective dimensionality condition. A generalization of these results to a wider class of effective dimensionality conditions and extension to kernel principal component regression can be found in \citet{Dicker17}.

For a variant of KCG \citet{Blan10b} obtained probabilistic convergence rates for independent identically distributed data. The pointed explicitly out that their approach and results are not directly applicable to KPLS.

We study of the convergence of the kernel partial least squares estimator to the true regression function when the algorithm is stopped early. 
Similar to \citet{Blan10b} we derive explicit probabilistic convergence rates. In contrast to previously cited works on kernel regression our input data are not independent and identically distributed but rather stationary time series. We derive probabilistic convergence results that can be applied for arbitrary temporal dependence structures, given that certain concentration inequalities for these data hold.
The derived convergence rates depend not only on the complexity of the target function and of the data mapped into the kernel space, but also on the persistence of the dependence in the data. In the stationary setting we prove that the short range dependence still leads to optimal rates, but if the dependence is more persistent, the rates become slower.

\section{Kernel Partial Least Squares}
\label{sec:problem}
Consider the non-parametric regression problem  
\begin{equation}
\label{eq:model}
	y_t = f^\ast(X_t) + \varepsilon_t,~~t \in \Z.
\end{equation}
Here $\{X_t\}_{t \in \Z}$ is a $d$-dimensional, $d \in \N$, stationary time series on a probability space $(\Omega,\calA,\Prob)$ and $\{\varepsilon_t\}_{t \in \Z}$ is an independent and identically distributed sequence of real valued random variables with expectation zero and variance $\sigma^2 > 0$ that is independent of $\{X_t\}_{t \in \Z}$. Let $X$ be a random vector that is independent of $\{X_t\}_{t \in \Z}$ and $\{\varepsilon_t\}_{t \in \Z}$ with the same distribution as $X_0$. The target function we seek to estimate is $f^\ast \in \LTwo$.

For the purpose of supervised learning assume that we have a training sample $\{(X_t,y_t)\}_{t=1}^n$ for some $n \in \N$. In the following we introduce some basic notation for the kernel based learning approach.

Define with $(\mathcal{H},\langle \cdot,\cdot\rangle_\calH)$ the RKHS of functions on $\R^d$ with reproducing kernel $k:\R^d \times \R^d \rightarrow \R$, i.e., it holds
\begin{equation}
\label{eq:rep.property}
	g(x) = \langle g, k(\cdot,x) \rangle_\calH, ~~x \in \R^d, g\in \calH.
\end{equation} 

The corresponding inner product and norm in $\calH$ is denoted by $\langle \cdot,\cdot \rangle_\calH$ and $\|\cdot\|_\calH$, respectively. We refer to \citet{Berlinet} for examples of Hilbert spaces and their reproducing kernels. In the following we deal with reproducing kernel Hilbert spaces which fulfill the following, rather standard, conditions:
\begin{enumerate}[label={(K\arabic*})]
\item \label{con:k1}
$\calH$ is separable,
\item \label{con:k2}
There exists a $\kappa>0$ such that $|k(x,y)| \leq \kappa$ for all $x,y \in \R^d$ and $k$ is measurable.
\end{enumerate}
Under \ref{con:k1} the Hilbert-Schmidt norm $\|\cdot\|_{\HS}$ for operators mapping from $\calH$ to $\calH$ is well defined.
If condition \ref{con:k2} holds, all functions in $\calH$ are bounded, see \citet{Berlinet}, chapter 2. 
The conditions are satisfied for a variety of popular kernels, e.g., Gaussian or triangular. 

The main principle of RKHS methods is the mapping of the data $X_t$ into $\calH$ via the feature maps $\phi_t = k(\cdot,X_t)$, $t=1,\dots,n$. This mapping can be done implicitly by using the kernel trick $\langle \phi_t,\phi_s \rangle_\calH = k(X_t,X_s)$ and thus only the $n \times n$ dimensional kernel matrix $K_n = n^{-1}[k(X_t,X_s)]_{t,s=1}^n$ is needed in the computations. Then the task for RKHS methods is to find coefficients $\alpha_1,\dots,\alpha_n$ such that $f_\alpha = \sum_{t=1}^n \alpha_t \phi_t$ is an adequate approximation of $f^\ast$ in $\calH$, measured in the $\LTwo$ norm $\|\cdot\|_2$.

There are a variety of different approaches to estimate the coefficients $\alpha_1,\dots,\alpha_n$, including kernel ridge regression, kernel principal component regression and, of course, kernel partial least squares. The latter method was introduced by \citet{Ros01} and is the focus of the current work.

It was shown by \citet{Krae07b} that the KPLS algorithm solves
\begin{equation}
\label{eq:kpls.optim}
	\alphaest_i = \arg\min\limits_{v \in \calK_i(K_n,y)} \|y - K_n v\|^2,~~i=1,\dots,n,
\end{equation}
with $y=(y_1,\dots,y_n)^\T$. Here $\calK_i(K_n,y) = \mathrm{span}\left\{
	y,K_n y, K_n^2y,\dots,K_n^{i-1}y
\right\}$, $i=1,\dots,n$, is the $i$th order Krylov space with respect to $K_n$ and $y$ and $\| \cdot\|$ denotes the Euclidean norm rescaled by $n^{-1}$. The dimension $i$ of the Krylov space is the regularization parameter for KPLS.

We will introduce several operators that will be crucial for our further analysis. Fist define two integral operators: the kernel integral operator $T^\ast:\LTwo \rightarrow \calH, g \mapsto\E\{k(\cdot,X) g(X)\}$ and the change of space operator $T:\calH \rightarrow \LTwo, g \mapsto g$, which is well defined if \ref{con:k2} holds.
It is easy to see that $T, T^\ast$ are adjoint, i.e., for $u \in \calH$ and $v \in \LTwo$ it holds $\langle T^\ast v, u \rangle_\calH = \langle v, T u \rangle_2$ with $\langle\cdot,\cdot \rangle_2$ being the inner product in $\LTwo$. 

The sample analogues of $T,T^\ast$ are $T_n:\calH \rightarrow \R^n, g \mapsto \{g(X_1),\dots,g(X_n)\}^\T$ and $T_n^\ast:\R^n \rightarrow \calH, (v_1,\dots,v_n)^\T \mapsto n^{-1} \sum_{t=1}^n v_t k(\cdot,X_t)$, respectively. Both operators are adjoint with respect to 
the rescaled Euclidean product $\langle u,v\rangle = n^{-1} u^\T v$, $u,v \in \R^d$

Finally, we define the sample kernel covariance operator $\KEst = T^\ast_n T_n:\calH \rightarrow \calH$ and the population kernel covariance operator $\KCov = T^\ast T:\calH \rightarrow \calH$. Note that it holds $K_n = T_n T_n^\ast$. Under \ref{con:k1} and \ref{con:k2} $\KCov$ is a self-adjoint compact operator with operator norm $\|\KCov\|_{\calL} \leq \kappa$, see \citet{Cap07}.

With this notation we can restate (\ref{eq:kpls.optim}) for the function $f_\alpha$ 
\begin{equation}
\label{eq:func.representation}
f_{\widehat{\alpha}_i} = \arg \min_{g \in \calK_i(\KEst,\bEst )} \|y- \{g(X_1),\dots,g(X_n)\}^\T\|^2=\arg \min_{g \in \calK_i(\KEst,\bEst )}\|y-T_ng\|^2.
\end{equation}
Hence, we are looking for functions that minimize the squared distance to $y$ constrained to a sequence of Krylov spaces. 

In the literature of ill-posed problems it is well known that without further conditions on the target function $f^\ast$ the convergence rate of the conjugate gradient algorithm can be arbitrarily slow, see \citet{Hanke}, chapter 3.2. One common a-priori assumption on the regression function $f^\ast$ is 
a source condition:
\begin{enumerate}[label={(S)}]
\item
\label{eq:source}
There exist $r \geq 0$, $R>0$ and $u \in \LTwo$ such that $f^\ast = (T T^\ast)^{r} u$ and $\|u\|_2 \leq R$.
\end{enumerate}

If $r \geq 1/2$, then the target function $f^\ast \in \LTwo$ coincides almost surely with a function $f \in \calH$ and we can write $f^\ast = T f$, see \citet{Cuc02}.
With this the kernel partial least squares estimator $f_{\widehat{\alpha}_i}$ estimates the correct target function, not only its best approximation in $\mathcal{H}$. This case is known as the inner case. 

The situation with $r<1/2$ is referred to as the outer case. Under additional assumptions, e.g., the availability of additional unlabeled data, it is still possible that an estimator of $f^\ast$ converges to the true target function in $\LTwo$ norm with optimal rates (with respect to the number $n$ of labeled data points). See \citet{Vit06} for a detailed description of this semi-supervised approach for kernel ridge regression in the independent and identically distributed case. We do not treat the case $r<1/2$ in this work.

A source conditions is often interpreted as an abstract smootheness condition. This can be seen as follows.
Let $\eta_1 \geq \eta_2 \geq \dots$ be the eigenvalues and $\psi_1,\psi_2,\dots$ the corresponding eigenfunctions of the compact operator $S$. 
Then it is easy to see that the source condition \ref{eq:source} is equivalent to $f = \sum_{j=1}^\infty b_j \psi_j \in \LTwo$ with $b_j$ such that $\sum_{j=1}^\infty \eta_j^{-2(r+1/2)} b_j^2 < \infty$. Hence, the higher $r$ is chosen the faster the sequence $\{b_j\}_{j=1}^\infty$ must converge to zero. Therefore, the sets of functions for which source conditions hold are nested, i.e., the larger $r$ is the smaller the corresponding set will be. The set with $r=1/2$ is the largest one and corresponds to a zero smoothness condition, i.e., $\sum_{j=1}^\infty \eta_j^{-2} b_j^2 < \infty$, which is equivalent to $f \in \calH$. For more details we refer to \citet{Dicker17}.

\section{Consistency of Kernel Partial Least Squares}
\label{sec:kpls.convergence}
The KCG algorithm as described by \citet{Blan10b} is consistent when stopped early and convergence rates can be obtained when a source condition \ref{eq:source} holds. Here we will proof the same property for KPLS. Early stopping in this context means that we stop the algorithm at some $a=a(n) \leq n$ and consider the estimator $f_{\alphaest_a}$ for $f^\ast$.

The difference between KCG and KPLS is the norm which is optimized. The kernel conjugate gradient algorithm studied in \citet{Blan10b} estimates the coefficients $\alpha \in \R^n$ of $f_\alpha$ via $\alphaest_i^{CG} = \arg\min_{v \in \calK_i(K_n,y)} \langle y- K_n v, K_n(y- K_n v)\rangle$.
It is easy to see that this optimization problem can be rewritten for the function $f_\alpha$ as  
\[
\min_{g \in \calK_i(\KEst,\bEst )} \|T_n^\ast y-S_ng\|_\calH^2=	\min_{g \in \calK_i(\KEst,\bEst )} \|T_n^\ast\left(y-T_ng\right)\|_\calH^2,
\] 
compared to (\ref{eq:func.representation}) for KPLS. Thus, KCG obtains the least squares approximation $g$ in the $\calH$-norm for the normal equation $\bEst = T^\ast _nT_n g$ and KPLS finds a function that minimizes the residual sum of squares. In both methods the solutions are restricted to functions $g \in \calK_i(\KEst,\bEst )$.

An advantage of the kernel conjugate gradient estimator is that concentration inequalities can be established for both $\bEst$ and $\KEst$ and applied directly as the optimization function contains both quantities. The stopping index for the regularization can be chosen by a discrepancy principle as $\aopt = \min\{1\leq i \leq n: \|\KEst f_{\alphaest_i^{CG}} - \bEst\| \leq \Lambda_n\}$ with $\Lambda_n$ being a threshold sequence that goes to zero as $n$ increases.

On the other hand, the function to be optimized for KPLS contains only $y$ and $T_n g = \{g(X_1),\dots,g(X_n)\}^\T$ for which statistical properties are not readily available. Thus, we need to find a way to apply the concentration inequalities for $\bEst$ and $\KEst$ to this slightly different problem. This leads to complications in the proof of consistency and a rather different and more technical stopping rule for choosing the optimal regularization parameter $a^\ast$ is used, as can be seen in Theorem \ref{th:kpls}. This stopping rule has its origin in \citet{Hanke}.\\

 In the following $\|\cdot\|_{\cal{L}}$ denotes the  operator norm and $\|\cdot\|_{HS}$ is the Hilbert-Schmidt norm.

\begin{theorem}
\label{th:kpls}
Assume that conditions \ref{con:k1}, \ref{con:k2}, \ref{eq:source} hold with $r \geq 3/2$ and there are constants $C_\delta(\nu),C_\epsilon(\nu)>0$ and a sequence $\{\gamma_n\}_{n \in \N} \subset [0,\infty)$, $\gamma_n \rightarrow 0$, such that we have for $\nu \in (0,1]$
\begin{align*}
	\Prob\left(
	\|\KEst-\KCov\|_{\calL} \leq C_\delta(\nu) \gamma_n
	\right) &\geq 1 - \nu/2,\\
	\Prob\left(
	\|\bEst - \KCov f\|_\calH \leq C_\epsilon(\nu) \gamma_n
	\right) &\geq 1- \nu/2.
\end{align*}
Define the stopping index $a^\ast$ by
\begin{equation}
\label{eq:stopping}
a^\ast = \min\left\{1 \leq a \leq n: \sum_{i=0}^a \|\KEst f_{\alphaest_i} - \bEst\|^{-2}_\calH \geq (C \gamma_n)^{-2}
	\right\},
\end{equation}
with $C = C_\epsilon(\nu) + \kappa^{r-1/2}(r+1/2) R \{1 + C_\delta(\nu)\}$.

Then it holds with probability at least $1-\nu$ that
\begin{align*}
	\|f_{\alphaest_{a^\ast}} - f^\ast\|_2 &= O\left\{\gamma_n^{2r/(2r+1)}\right\},\\
	\|f_{\alphaest_{a^\ast}} - f\|_\calH &= O\left\{\gamma_n^{(2r-1)/(2r+1)}\right\},
\end{align*}
with $f^\ast = T f$.
\end{theorem}

It can be shown that the stopping rule (\ref{eq:stopping}) always determines a finite index, i.e., the set the minimum is taken over is not empty, see \citet{Hanke}, chapter 4.3.

The theorem yields two convergence results, one in the $\calH$-norm and one in the $\LTwo$-norm. It holds that $\|v\|_2 = \|S^{1/2}v\|_\calH$. These are the endpoints of a continuum of norms $\|v\|_\beta = \|S^\beta v\|_\calH$, $\beta \in [0,1/2]$ that were considered in \citet{Nemirovskii86} for the derivation of convergence rates for KCG algorithms in a deterministic setting. 

The convergence rate of the kernel partial least squares estimator depends crucially on the sequence $\gamma_n$ and the source parameter $r$. If $\gamma_n = O(n^{-1/2})$, this yields the same convergence rate as Theorem 2.1 of \citet{Blan10b} for kernel conjugate gradient or \citet{Vito05} for kernel ridge regression with independent and identically distributed data. For stationary Gaussian time series we will derive concentration inequalities in the next section and obtain convergence rates depending on the source parameter $r$ and the range of dependence. Note that Theorem \ref{th:kpls} is rather general and it can be applied to any kind of dependence structure, as long as the necessary concentration inequalities can be established.

The next theorem  derives faster convergence rates under assumptions on the effective dimensionality of operator $S$, which is defined as $\ed = \mathrm{tr}\{(S+\lambda)^{-1} S\}$. The concept of effective dimensionality was introduced in \citet{Zho02} to get sharp error bounds for general learning problems considered there. If $\calH$ is a finite dimensional space it was shown in \citet{Zho02} that $d_\lambda \leq \mathrm{dim}(\calH)$. For infinite dimensional spaces it describes the complexity of the interactions between data and reproducing kernel. 

If $d_\lambda = O(\lambda^{-s})$ for some $s \in (0,1]$, \citet{Cap07} showed that the order optimal convergence rates $n^{-r/(2r+s)}$ are attained for KRR with independent and identically distributed data.

The effective dimensionality clearly depends on the behaviour of eigenvalues of $S$. If these converge sufficiently fast to zero, nearly parametric rates of convergence can be achieved for reproducing kernel Hilbert space methods, see, e.g., \citet{Dicker17}. In particular, the behaviour of $d_\lambda$ around zero is of interest, since it determines how ill-conditioned the operator $(S+\lambda)^{-1}$ becomes. In the following theorem we set $\lambda = \lambda_n$ for a sequence $\{\lambda_n\}_{n\in \N} \subset (0,\infty)$ that converges to zero. 

\begin{theorem}
\label{th:kpls2}
Assume that conditions \ref{con:k1}, \ref{con:k2}, \ref{eq:source} hold with $r \geq 1/2$ and that the effective dimensionality $\ed$ is known. Additionally, there are constants $C_\delta(\nu),C_\epsilon(\nu),C_\psi
>0$ and a sequence $\{\gamma_n\}_{n \in \N} \subset [0,\infty)$, $\gamma_n \rightarrow 0$, such that for $\nu \in (0,1]$ and $n$ sufficiently large
\begin{align*}
	\Prob\left\{
	\|\KEst-\KCov\|_{\calL} \leq C_\delta(\nu) \gamma_n
	\right\} &\geq 1 - \nu/3,\\
	\Prob\left\{
	\|(S+\lambda_n)^{-1/2}(\bEst - \KCov f)\|_\calH \leq C_\epsilon(\nu) \sqrt{d_{\lambda_n}}\gamma_n
	\right\} &\geq 1- \nu/3,\\
	\Prob\left\{
	\|(S+\lambda_n)^{1/2}(S_n+\lambda_n)^{-1/2}\|_{\calL} \leq C_\psi
	\right\} &\geq 1 - \nu/3,
\end{align*}
Here $\{\lambda_n\}_{n \in \N} \subset (0,\infty)$ is a sequence converging to zero such that for $n$ large enough 
\begin{equation}
	\label{eq:lambda.inequ}
	\gamma_n \leq \lambda_n^{r-1/2}.
\end{equation}

Take $\zeta_n = \max\{\sqrt{\lambda_n d_{\lambda_n}} \gamma_n, \lambda_n^{r+1/2}\}$ 
Define the stopping index $a^\ast$ by
\begin{equation}
\label{eq:stopping2}
a^\ast = \min\left\{1 \leq a \leq n: \sum_{i=0}^a \|\KEst f_{\alphaest_i} - \bEst\|^{-2}_\calH \geq (C \zeta_n)^{-2}\right\},
\end{equation}
with $C=4 R \max\{1, C_\psi^2,(r-1/2)\kappa^{r-3/2} C_\delta(\nu),2^{-1/2}R^{-1} C_\psi C_\epsilon(\nu)\}$.

Then it holds with probability at least $1-\nu$ that
\begin{align*}
	\|f_{\alphaest_{a^\ast}} - f^\ast\|_2 &= O\left\{\lambda_n^{-1/2}\zeta_n\right\},\\
\|f_{\alphaest_{a^\ast}} - f\|_\calH &= O\left\{\lambda_n^{-1}\zeta_n\right\},
\end{align*}
with $f^\ast = T f$.
\end{theorem}
The condition (\ref{eq:lambda.inequ}) holds trivially for $r=1/2$ as $\gamma_n$ converges to zero. For $r >1/2$ the sequence $\lambda_n$ must not converge to zero arbitrarily fast. 

In its general form Theorem \ref{th:kpls2} does not give immediate insight in the probabilistic convergence rates of the kernel partial least squares estimator. Therefore, we state two corollaries, where the function $d_\lambda$ is specified. In both corollaries we explicitly state the choice of the sequence $\lambda_n$ that yield the corresponding rates.

\begin{corollary}
\label{cor:pol.ed}
Assume that there exists $s \in (0,1]$ such that $
	d_\lambda = O(\lambda^{-s})
$ for $\lambda \rightarrow 0$.
Then under conditions of Theorem \ref{th:kpls2} with $\lambda_n = \gamma_n^{2/(2r + s)}$  it holds with probability at least $1-\nu$ that
\begin{align*}
	\|f_{\alphaest_{a^\ast}} - f^\ast\|_2 &= O\left\{\gamma_n^{2r/(2r+s)}\right\}.
\end{align*}
\end{corollary}
Polynomial decay of the effective dimensionality $d_\lambda = \mathrm{tr}\{(S+\lambda)^{-1}S\}$ occurs if the eigenvalues of $S$ also decay polynomially fast, that is, $\mu_i = c_s i^{-1/s}$ for $s \in (0,1]$, since in this case $d_\lambda = \sum\limits_{i=1}^\infty \{1+\lambda/c_s i^{1/s}\}^{-1} = O(\lambda^{-s})$. This holds, for example, for the Sobolev kernel $k(x,y) = \min(x,y)$, $x,y \in [0,1]$ and data that are uniformly distributed on $[0,1]$, see \citet{Ras14}.

If $\gamma_n = n^{-1/2}$, then the KPLS estimator converges in the $\LTwo$-norm with a rate of $n^{-r/(2r+s)}$. This rate is shown to be optimal in \citet{Cap07} for KRR with independent identically distributed data.

Note that the rate obtained in Theorem \ref{th:kpls} corresponds to $\gamma_n^{-2r/(2r+s)}$ with $s=1$, i.e., the worst case rate with respect to the parameter $s \in (0,1]$.

In the next corollary to Theorem \ref{th:kpls2} we assume that the effective dimensionality behaves in a logarithmic fashion.

\begin{corollary}
\label{cor:log.ed}
Let $d_\lambda = O\{\log(1+a/\lambda)\}$ for $\lambda \rightarrow 0$ and  $a>0$. Then under the conditions of Theorem \ref{th:kpls2} with $\lambda_n = \gamma_n^2 \log\{ \gamma_n^{-2}\}$ and $r=1/2$ it holds with probability at least $1-\nu$ that
\begin{align*}
	\|f_{\alphaest_{a^\ast}} - f^\ast\|_2 &= O\left\{\gamma_n \log(1/2 \gamma_n^{-2})\right\}.
\end{align*}
\end{corollary}
The effective dimensionality takes the special form considered in this corollary, for example, when the eigenvalues of $S$ decay exponentially fast. This holds, for example, if the data are Gaussian and the Gaussian kernel is used, see Section \ref{sec:gauss.kern}.
If $\gamma_n =O(n^{-1/2})$, then the convergence rate is of order $O\{n^{-1}\log(n)\}$, which are nearly parametric. It is noteworthy that the source condition only impacts the choice of the sequence $\lambda_n$, not the convergence rates of the estimator in the $\LTwo$-norm. Therefore, we stated the corollary for $r=1/2$, which is a minimal smoothness condition on $f^\ast$, i.e., that $f^\ast = T f$ almost surely for an $f \in \calH$.

The rates obtained in Corollaries \ref{cor:pol.ed} and \ref{cor:log.ed} were derived in \citet{Dicker17} for kernel ridge regression and kernel principal component regression under the assumption of independent and identically distributed data.

\section{Concentration Inequalities for Gaussian Time Series}
\label{sec:concentration}
Crucial assumptions of Theorem \ref{th:kpls} and \ref{th:kpls2} are the concentration inequalities
for $\KEst$ and $\bEst$ and convergence of the sequence $\{\gamma_n\}_{n \in \N}$. Here we establish such inequalities in a Gaussian setting for stationary time series. At the end of this section we will state explicit convergence rates for $f_{\widehat{\alpha}_\aopt}$ that depend not only on the source parameter $r \geq 1/2$ and the effective dimensionality $\ed$, but also on the persistence of the  dependence in the data.

The Gaussian setting is summarized in the following assumptions
\begin{enumerate}[label={(D\arabic*})]
\item
\label{D1}
$(X_h,X_0)^\T \sim \mathcal{N}_{2d}(0,\NCov_h)$, $h=1,\dots,n-1$, with 
\[
      \NCov_h = \left[
\begin{matrix}
	\tau_0 & \tau_h\\
	\tau_h & \tau_0
\end{matrix}
\right] \otimes \NCov.
\]
Here $\NCov \in \R^{d \times d}$ and $V=[\tau_{|i-j|}]_{i,j=1}^{n} \in \R^{n \times n}$ are positive definite, symmetric matrices and $\otimes$ denotes the Kronecker product between matrices. Furthermore $X_0 \sim \mathcal{N}_d(0,\tau_0 \Sigma)$.
\item
\label{D2}
For the autocorrelation function $\rho_h = \tau^{-1}_0\tau_h$ there exists a $q>0$ such that $|\rho_h| \leq  (h+1)^{-q}$ for $h =0,\dots,n-1$.
\end{enumerate}
Condition \ref{D1} is a separability condition for the covariance matrices $\Sigma_h$, $h = 0,\dots,n-1$. Due to \ref{D1} the effects (on the covariance) over time and between the different variables can be treated separately. Under condition \ref{D2} it is easy to see that from $q>1$ follows the absolute summability of the autocorrelation function $\rho$ and thus $\{X_t\}_{t \in \Z}$ is a short memory process. Stationary short memory processes keep many of the properties of independent and identically distributed data, see, e.g., \citet{bBrock}.

On the other hand $q \in (0,1]$ yields a long memory process, see, e.g., Definition 3.1.2 in \citet{Giraitis}. Examples of long memory processes are the fractional Gaussian noise with an autocorrelation function that behaves like $(h+1)^{-2(1-H)}$, with $H \in [0,1)$ being the Hurst coefficient. Stationary long memory processes exhibit dependencies between observations that are more persistent and many statistical results that hold for independent and identically distributed data turn out to be false, see \citet{Samoro} for details.

The next theorem gives concentration inequalities for both estimators $\KEst$ and $\bEst$ in a Gaussian setting with convergence rates depending on the parameter $q>0$. These inequalities are the ones needed in Theorem \ref{th:kpls} and Theorem \ref{th:kpls2}. Recall that $\ed = \mathrm{tr}\{(S+\lambda)^{-1}S\}$ denotes the effective dimensionality of $S$.

\begin{theorem}
\label{th:conc.equality}
(i) Define $\diff\mu_h(x,y) = \diff\Prob^{X_h,X_0}(x,y) - \diff \Prob^{X_0}(x)\diff \Prob^{X_0}(y)$.
Under Assumptions \ref{con:k1} and \ref{con:k2} it holds for $\nu \in (0,1]$ with probability at least $1-\nu$ that
\begin{align*}
\|\KEst - \KCov\|^2_{\calL} &\leq
\frac{2 \nu^{-1}}{n^2}\sum\limits_{h=1}^{n-1} (n-h) \int\limits_{\R^{2d}} k^2(x,y) \diff\mu_h(x,y) + \frac{\nu^{-1}}{n} \left\{
	\E k^2(X_0,X_0) - \|\KCov\|^2_{\HS}
\right\},\\
\|T_n^\ast y - \KCov f\|^2_\calH &\leq
\frac{2\nu^{-1}}{n^2}\sum\limits_{h=1}^{n-1} (n-h) \int\limits_{\R^{2d}} k(x,y)f(x)f(y)\diff\mu_h(x,y)\\
&+ \frac{\nu^{-1}}{n} \left[
	\E \left\{k(X_0,X_0) f^2(X_0)\right\} - \| \KCov f\|^2_\calH + \sigma^2 \E\{ k(X_0,X_0)\}
\right].
\end{align*}

(ii) Assume that additionally to \ref{con:k1}, \ref{con:k2} also \ref{D1}, \ref{D2} for $q >0$ are fulfilled. Denote $M=\sup_{x \in \R^d} |f(x)|$.

Then there exists a constant $C(q)>0$ such that 
\begin{align*}
\|\KEst - \KCov\|_\calL
	&\leq 
	 \nu^{-1/2} \{\gamma_n^2(q) \kappa C_\gamma	
		+ n^{-1}(\kappa^2-\|S\|_\HS^2)\}^{1/2},\\
 \|T_n^\ast y - \KCov f\|_\calH
	&\leq 
	\nu^{-1/2} \left[\gamma_n^2(q) M C_\gamma 
		+ n^{-1}\left\{
		\kappa (M + \sigma^2) - \| \KCov f\|^2_\calH	\right\}\right]^{1/2}
		,
\end{align*}
for $C_\gamma = C(q)\{(2\pi)^d \mathrm{det}(\Sigma)\}^{-1/2} \kappa d^{1/2}(1-4^{-q})^{-1/4(d+2)}$. The function $\gamma_n(q)$, $q>0$, is defined as
\[
	\gamma_n(q) = \left\{
				\begin{array}{clc}
			 n^{-1/2} &,& q>1\\
		n^{-1/2} \log(1/2n)	&,& q=1\\
			n^{-q/2} 
			&,& q \in (0,1).
		\end{array}
	\right.
\]

(iii) Let \ref{con:k1}, \ref{con:k2} and \ref{eq:source} hold. Let $\gamma_n(q)$ be the function as defined in (ii). Then there exists a constant $\tilde{C}_\epsilon>0$ such that it holds with probability at least $1-\nu$ for $\lambda>0$ that
\[
 \|(S+\lambda)^{-1/2}(T_n^\ast y - \KEst f\|_\calH \leq  \nu^{-1/2}  \tilde{C}_\epsilon \sigma \sqrt{\ed}\gamma_n(q).
\] 

(iv) Let \ref{con:k1}, \ref{con:k2}, \ref{eq:source}, \ref{D1} and \ref{D2} hold. Let $\lambda_n^{-1/2}d_{\lambda_n}^{1/2} \gamma_n(q) \rightarrow 0$ for a sequence $\lambda_n \rightarrow 0$ and $\gamma_n(q)$ the function defined in (ii). Then there exists an $n_0 = n_0(\nu,q) \in \N$ such that with probability at least $1-\nu$ we have for all $n \geq n_0$
\[\|(S+\lambda_n)^{1/2}(\KEst+\lambda_n)^{-1/2}\|_\calL \leq \sqrt{2}.
\]
\end{theorem}

The first part of the theorem is general and can be used to derive concentration inequalities not only in the Gaussian setting and is of interest in itself. The convergence rate is controlled by the sums appearing on the right hand side.
 If these sums are of $O(n)$ then the mean squared error of both $\KEst$ and $\bEst$ will converge to zero with a rate of $n^{-1}$. On the other hand, if the sums are of order $O(n^{2-q})$ for some $q\in (0,1)$, the mean squared errors will converge with the reduced rate $n^{-q}$.

The second part derives explicit concentration inequalities in the Gaussian setting described by \ref{D1} and \ref{D2} with rates depending on the range of the dependence measured by $q>0$. These inequalities appear in Theorem \ref{th:kpls}.

Parts (iii) and (iv) give the additional probabilistic bounds needed to apply Theorem \ref{th:kpls2}. The condition $\lambda_n^{-1/2} d^{1/2}_{\lambda_n} \gamma_n(q) \rightarrow 0$ in Theorem \ref{th:conc.equality} (iv) is fulfilled in the settings of Corollary \ref{cor:pol.ed} and Corollary \ref{cor:log.ed}.

 Theorem \ref{th:kpls}, Corollary \ref{cor:pol.ed}, Corollary \ref{cor:log.ed} and Theorem \ref{th:conc.equality} together imply
\begin{corollary}
\label{cor:convergence}
Let the conditions of Theorem \ref{th:kpls2} and \ref{D1}, \ref{D2} hold.

(i) Assume that there exists $s\in (0,1]$ such that $d_\lambda = O(\lambda^{-s})$ for $\lambda \rightarrow 0$. Then with probability at least $1-\nu$
	\[
		\|f_{\alphaest_{a^\ast}} - f^\ast\|_2 = \left\{
			\begin{array}{cc}
				O\{n^{-r/(2r+s)}\}, & q>1,\\
				O\{n^{-q r/(2r+s)}\}, & q \in (0,1).
			\end{array}
			\right.
	\]
If instead of conditions of Theorem \ref{th:kpls2}, conditions of Theorem \ref{th:kpls} are assumed, then the convergence rates  above have $s=1$. 
	
(ii) Assume that there exists $a>0$ such that $\ed = O\{\log(1+a/\lambda)\}$ for $\lambda \rightarrow 0$ and $r=1/2$. Then with probability at least $1-\nu$
	\[
		\|f_{\alphaest_{a^\ast}} - f^\ast\|_2 = \left\{
			\begin{array}{cc}
				O\{n^{-1/2}\log(1/2 n)\}, & q>1,\\
				O\{n^{-q/2} \log(1/2 n^q)\}, & q \in (0,1).
			\end{array}
			\right.
	\]
\end{corollary}
Hence, for $q>1$ the kernel partial least squares algorithm achieves the same rates as if the data were independent and identically distributed. For $q \in (0,1)$ the convergence rates become substantially slower, highlighting that dependence structures that persist over a long time can influence the convergence rates of the algorithm.

\section{Source condition and effective dimensionality for Gaussian kernels}
\label{sec:gauss.kern}
The source condition \ref{eq:source} and the effective dimensionality $d_\lambda$ are of great importance in the  convergence rates derived in previous sections. Here we investigate these conditions for the reproducing kernel Hilbert space corresponding to the Gaussian kernel $k(x,y) = \exp(-l\|x-y\|^2)$, $x,y \in \R^d$, $l>0$, for $d=1$. Hence, the space $\calH$ is the space of all analytic functions that decay exponentially fast, see \citet{Steinwart05anexplicit}.

We also impose the normality conditions \ref{D1} and \ref{D2} on $\{X_t\}_{t\in \Z}$, where now $\sigma^2_x =\Sigma \in \R$ due to $d=1$. The following proposition derives a more explicit representation for $f \in \calH$.

\begin{proposition}
\label{prop:source}
Assume that \ref{con:k1},\ref{con:k2} and \ref{eq:source} hold for $r \geq 1/2$. Let $d=1$, $X_0 \sim \mathcal{N}(0,\sigma^2_x), \sigma^2_x > 0$ and consider the Gaussian kernel $k(x,y) = \exp\{-l (x-y)^2\}$ for $x,y \in \R$, $l>0$. Then $f$ can be expressed for $\mu = r - 1/2 \in \N$
  via $f(x) = \sum_{i=1}^\infty c_i L_\mu(x,z_i)$ for fixed $\{z_i\}_{i=1}^\infty,\{c_i\}_{i=1}^\infty \subset \R$ such that $\sum_{i,j=1}^\infty c_i c_j k(z_i,z_j) \leq R^2$, $R>0$. Here we have for $x,z \in \R$ 
	\begin{align*}
	L_\mu(x,z) &=  \exp
		\left[
			-1/2 \left\{
				\frac{\det(\Lambda)(x^2+z^2)-2 l^{\mu+1} x z}{\det(\Lambda_{1:\mu})}
			\right\}
		\right],
		\end{align*}
	with $\Lambda \in \R^{(\mu+1)\times (\mu+1)}$ being a tridiagonal matrix with elements 
	\[
		\Lambda_{i,j} = \left\{
	\begin{array}{cll}
		\sigma^{-2}_x + 2 l &, & i=j<\mu+1\\
		l &, & i=j=\mu+1\\
		-l &, & |i-j| = 1\\
		0&, & else
	\end{array}\right.
	\]
	for $i,j=1,\dots,\mu+1$ and $\Lambda_{1:\mu}$ is the $\mu \times \mu$-dimensional sub-matrix of $\Lambda$ including the fist $\mu$ columns and rows.
	
	Conversely any function $f^\ast = T f$ with $f$ of the above form fulfills a source condition $\ref{eq:source}$ with $r = \mu + 1/2$, $\mu \in \N$.
\end{proposition}
Hence if we fix an $r \geq 1/2$ with $r-1/2 \in \N$ this theorem gives us a way to construct functions $f \in \calH$ with $f^\ast = Tf$ that fulfill \ref{eq:source}.

The next proposition derives the effective dimensionality $\ed$ in this setting:
\begin{proposition}
\label{prop:ed}
Let $d=1$, $X_0 \sim \mathcal{N}(0,\sigma^2_x)$ for some $\sigma^2_x>0$ and consider the Gaussian kernel $k(x,y) = \exp\{-l(x-y)^2\}$, $x,y \in \R$, $l>0$. 

Then there is a constant $D>0$ such that it holds for any $\lambda \in (0,1]$
\[
	d_\lambda = \mathrm{tr}\{(S+\lambda)^{-1} S\} \leq 	D\log(1+a/\lambda),
\]
with $a =\sqrt{2}(1+\beta+\sqrt{1+\beta})^{-1/2} $, $\beta = 4 l \sigma^2_x$.
\end{proposition}

With the latter result Corollary \ref{cor:log.ed} is applicable and we expect convergence rates for the kernel partial least squares algorithm of order $O\{\gamma_n \log(1/2 \gamma_n^{-2})\}$ for a sequence $\{\gamma_n\}_n$ as in Theorem \ref{th:kpls2}.

\section{Simulations}
\label{sec:simulations}
To validate the theoretical results of the previous sections we conducted a simulation study. The reproducing kernel Hilbert space is chosen to correspond to the Gaussian kernel $k(x,y) = \exp(-l\|x-y\|^2)$, $x,y \in \R^d$, $l=2$, for $d=1$.

\begin{figure}
   \begin{center}
   	\includegraphics[width=.5\linewidth]{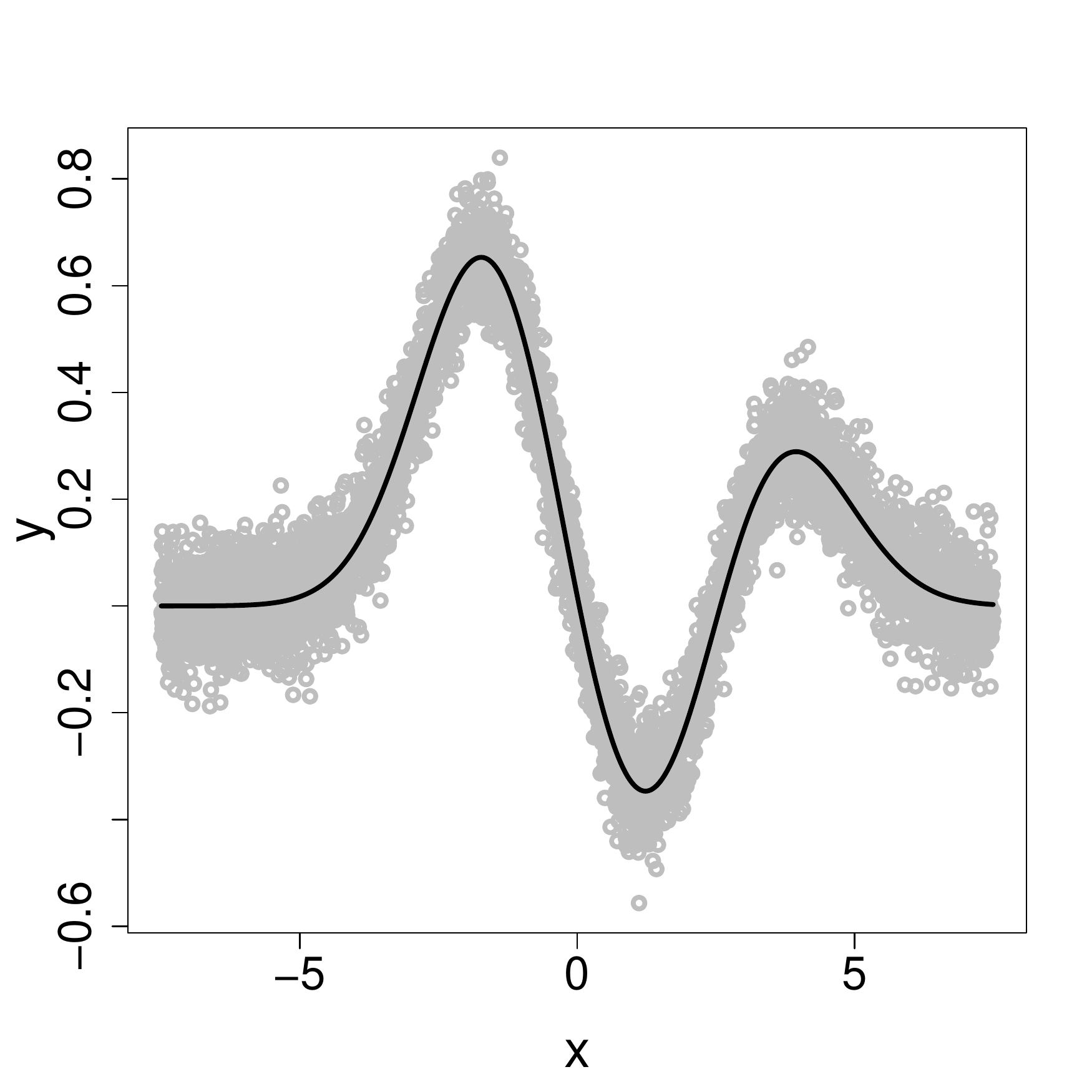}
   	\caption{\label{fig:func}
		The function $f$ evaluated on $[-7.5,7.5]$ (black) and one realisation of the noisy data $y = f(x) + \varepsilon$ (grey).   		
   		   }
   \end{center}
\end{figure}

The source parameter is taken $r=4.5$ and we consider the function 
\[
	f(x) = 4.37^{-1}\{3 {L}_4(x,-4) - 2 {L}_4(x,3)+ 1.5 {L}_4(x,9)\}, ~~ x \in \R.
	\] 
	The normalization constant is chosen such that $f$ takes values in $[-0.35,0.65]$ and $L_4$ is the exponential function given in Proposition \ref{prop:source}. The function $f$ is shown in Figure \ref{fig:func}.

In condition \ref{D1} we set $\sigma^2_x =\NCov = 4$. For the matrix $V^2 = [\tau_{|i-j|}]_{i,j=1}^n \in \R^{n \times n}$ we choose three different structures  for $n\in\{200,400,1000\}$. In the first setting $\tau_h = \mathbb{I}(h=0)$, which corresponds to independent data. The second setting with $\tau_h = 0.9^{-h}$ implies an autoregressive process of order one. Finally, the third setting with $\tau_h = (1+h)^{-q}$, $q=1/4$, $h =0,\dots,n-1$ leads to the long range dependent case. 

In a Monte Carlo simulation with $M=1000$ repetitions the time series $\{X_t^{(j)}\}_{t=1}^n$ are generated via $X^{(j)} = V N^{(j)}$ with $N^{(j)} \sim \mathcal{N}_n(0,\sigma^2 I_n)$, $j=1,\dots,M$, where $I_n$ is the $n \times n$-dimensional identity matrix. 

The residuals $\varepsilon_1^{(j)},\dots,\varepsilon_n^{(j)}$ are generated as independent standard normally distributed random variables and independent of $\{X_t^{(j)}\}_{t=1}^n$ . The response is defined as $y_t^{(j)} = f(X_t^{(j)}) + \eta\, \varepsilon_t^{(j)}$, $t=1,\dots,n$, $j=1,\dots,M$, with $\eta = 1/16$.

The kernel partial least squares and kernel conjugate gradient algorithms are run for each sample $\{(X_t^{(j)},y_t^{(j)})^\T\}_{t=1}^n$, $j=1,\dots,M$, with a maximum of $40$ iteration steps. We denote the estimated coefficients with $\alphaest_1^{(j,m)},\dots,\alphaest_{40}^{(j,m)}$, $j=1,\dots,M$, with $m=CG$ meaning that the kernel conjugate gradient algorithm was employed and $m=PLS$ that kernel partial least squares was used to estimate $\alpha_1,\dots,\alpha_n$. 

The squared error in the $\LTwo$-norm is calculated via
\[
	\widehat{e}_{n,\tau}^{(j,m)} =  \min\limits_{a=1,\dots,40} \left[ \frac{1}{\sqrt{2\pi\sigma_x^2}}\int\limits_{-\infty}^\infty \left\{f_{\widehat{\alpha}_a^{(j,m)}}(x) - f(x) \right\}^2 \exp
	\left(-\frac{1}{2\sigma_x^2} x^2 \right) \diff x \right],
\]
for $j=1,\dots,M$, $n = 200,400,\dots,1000$ and $m \in \{CG, PLS\}$.

\begin{figure}
   \begin{center}
	 	\begin{minipage}{.33\textwidth}
	  		\centering
	 		\includegraphics[width=\ssize]{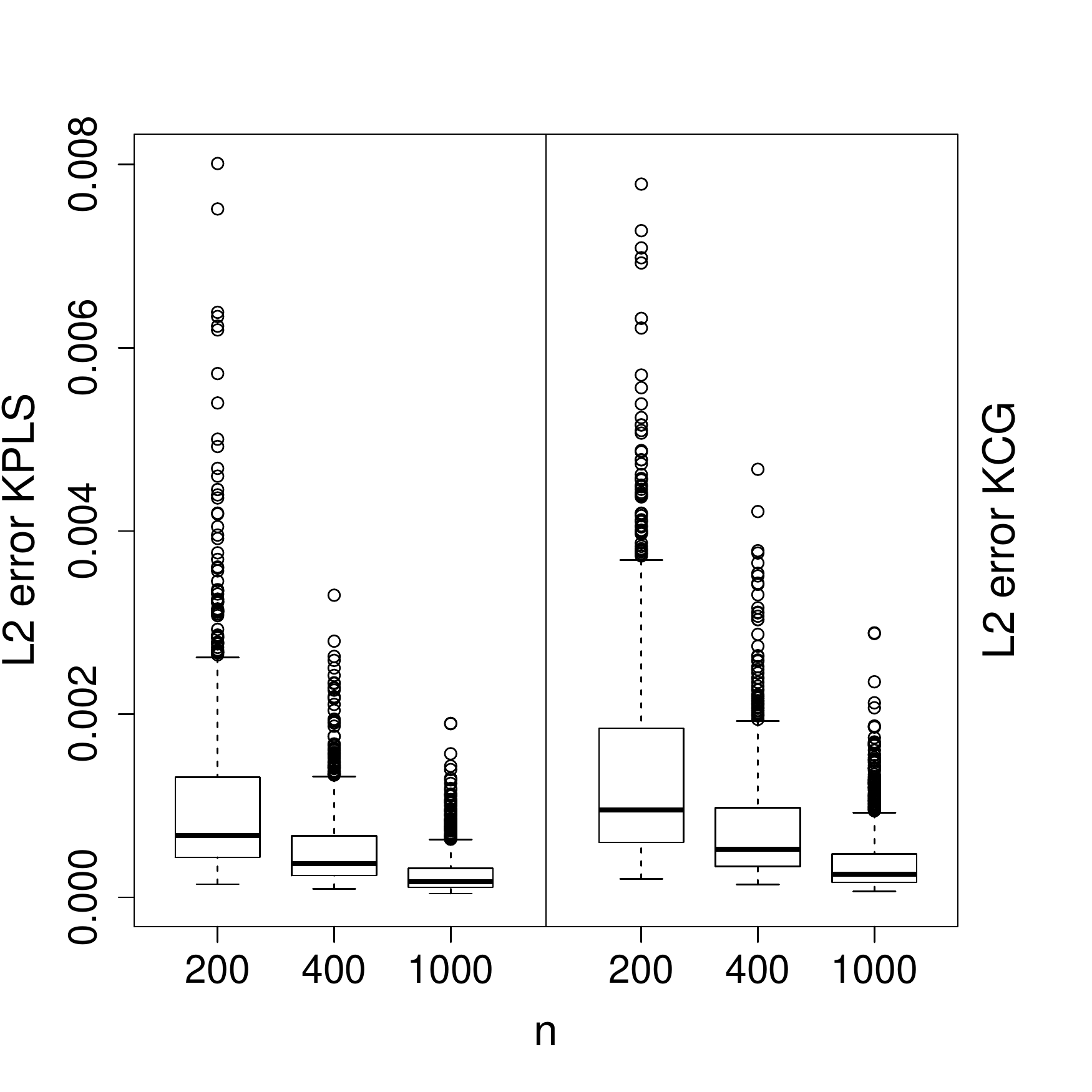}
		\end{minipage}%
		\begin{minipage}{.33\textwidth}
	  		\centering
	 		\includegraphics[width=\ssize]{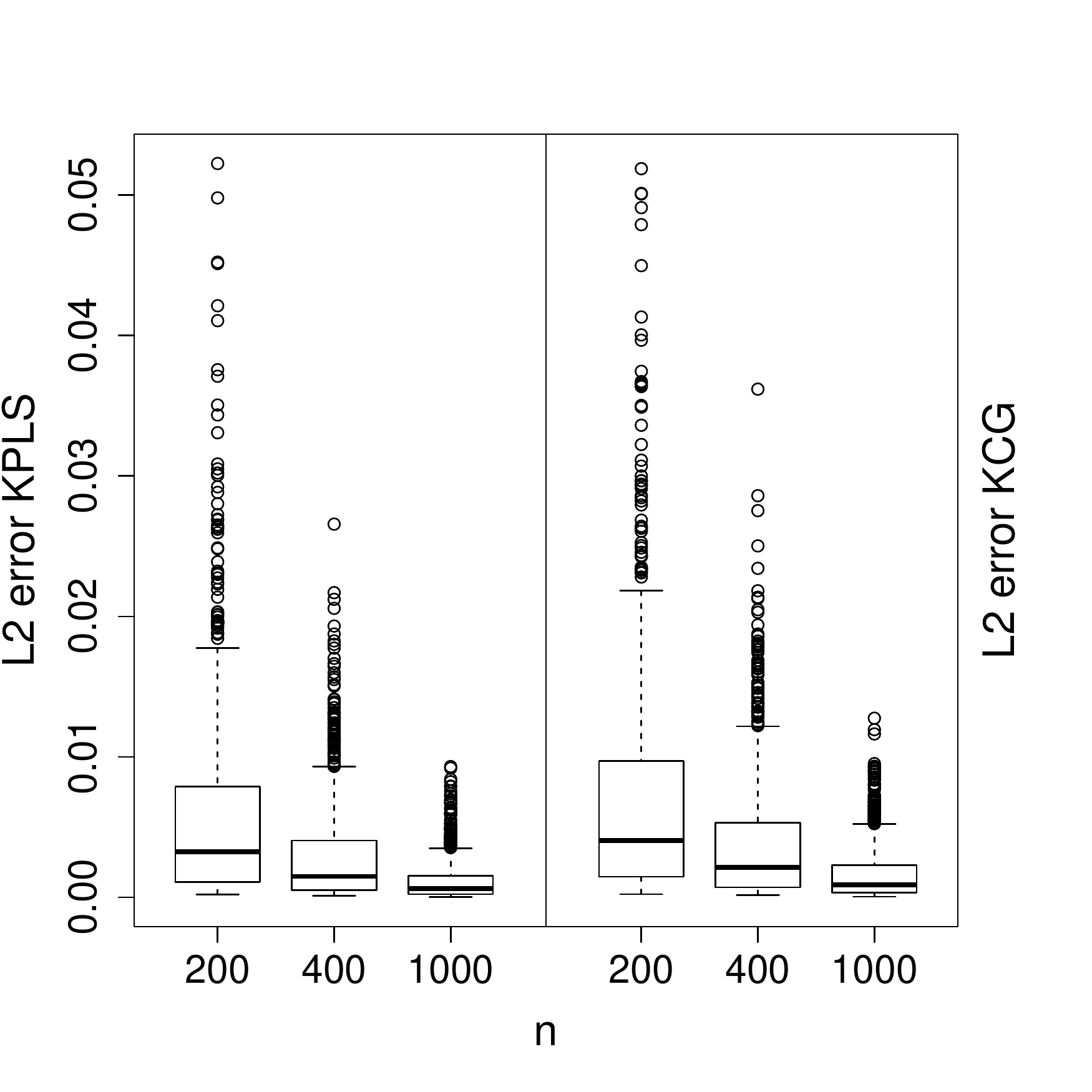}
		\end{minipage}%
		\begin{minipage}{.33\textwidth}
	  		\centering
	 		\includegraphics[width=\ssize]{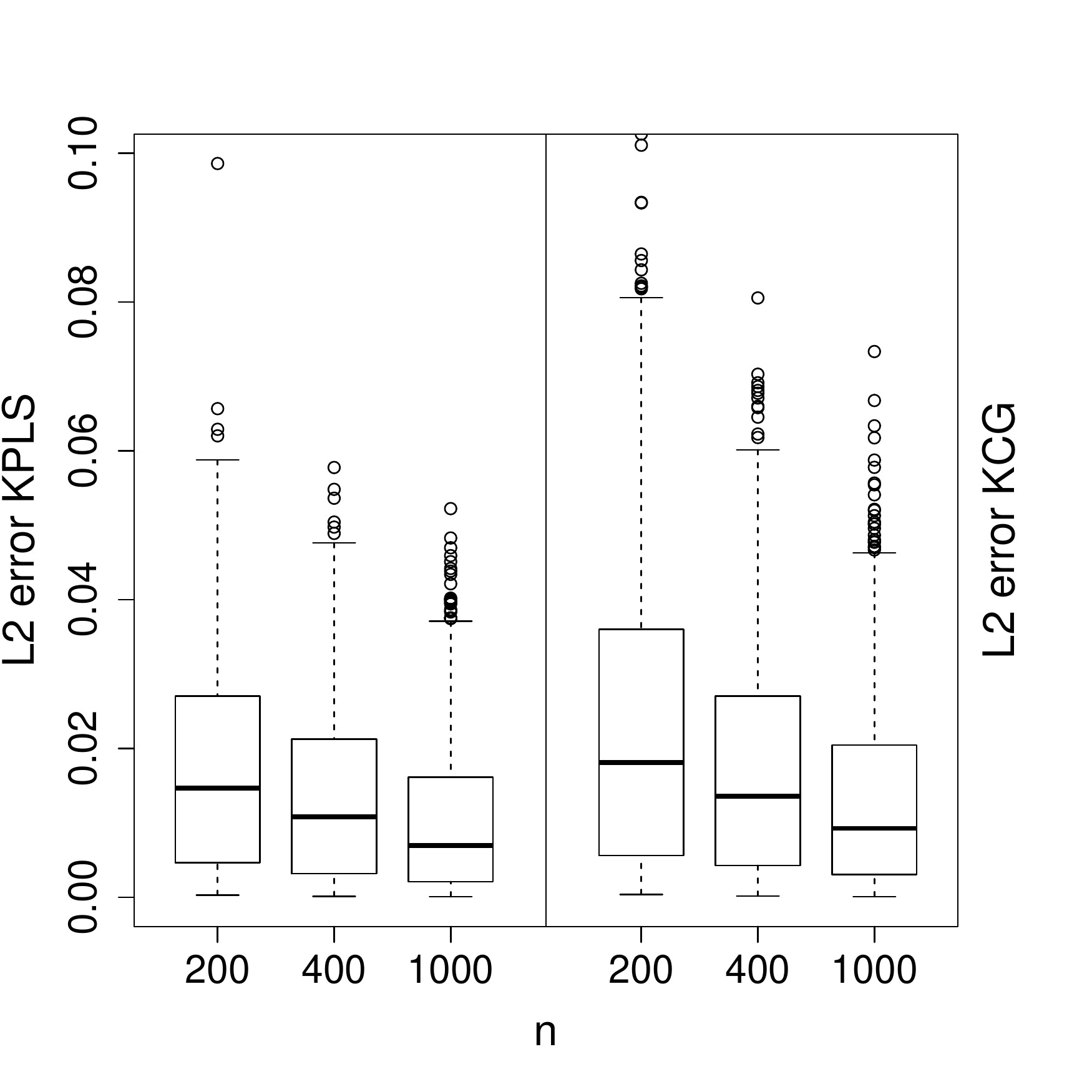}
		\end{minipage}
		\caption{\label{fig:box}
			Boxplots of the $\LTwo$-errors $\{\widehat{e}_{n,\tau}^{(j,m)}\}_{j=1}^{M}$ of kernel partial least squares (left side of each panel) and kernel conjugate gradient (right side of each panel) for different autocovariance functions $\tau$ and $n = 200,400,1000$. On the left is $\tau_h = \mathbb{I}(h=0)$, in the middle $\tau_h = 0.9^{-h}$ and on 
			the right $\tau_h = (h+1)^{-1/4}$.
		}
	\end{center}
\end{figure}

The results of the Monte-Carlo simulations are depicted in the boxplots of Figure \ref{fig:box}.
For kernel partial least squares (left panels) one observes that independent and autoregressive dependent data have roughly the same convergence rates, although the latter have a somewhat higher error. In contrast, the long range dependent data show slower convergence with the larger interquartile range, supporting the theoretical results of Corollary \ref{cor:convergence}.

The $\LTwo$-error of kernel conjugate gradient estimators is generally slightly higher than that of kernel partial least squares. Nonetheless, both of them have a similar behaviour.

\begin{figure}
   \begin{center}
   	\begin{minipage}{.45\textwidth}
	\center   	
   	\includegraphics[width=\ssize]{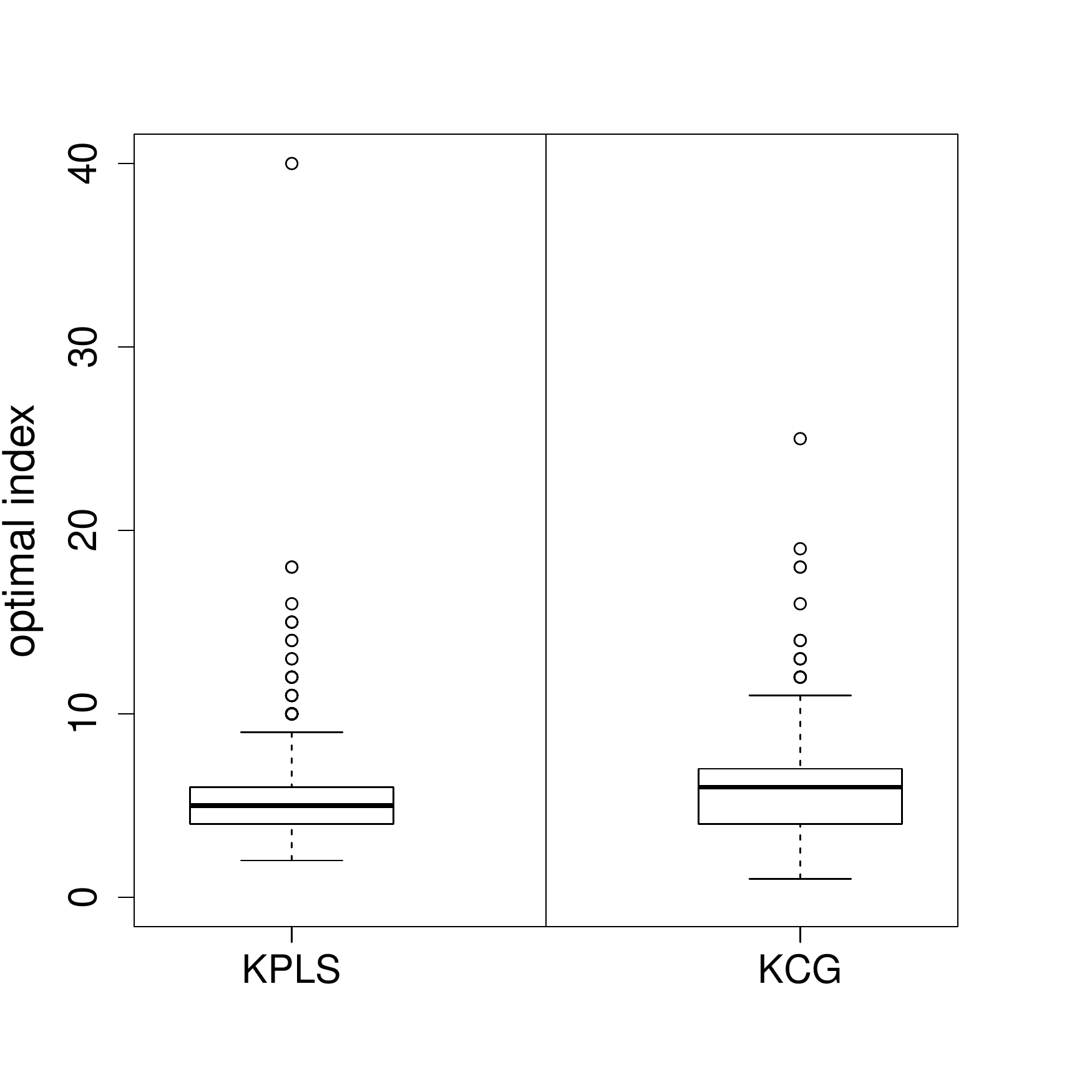}
   	\end{minipage}
   	\begin{minipage}{.45\textwidth}
	\center   	
   	\includegraphics[width=\ssize]{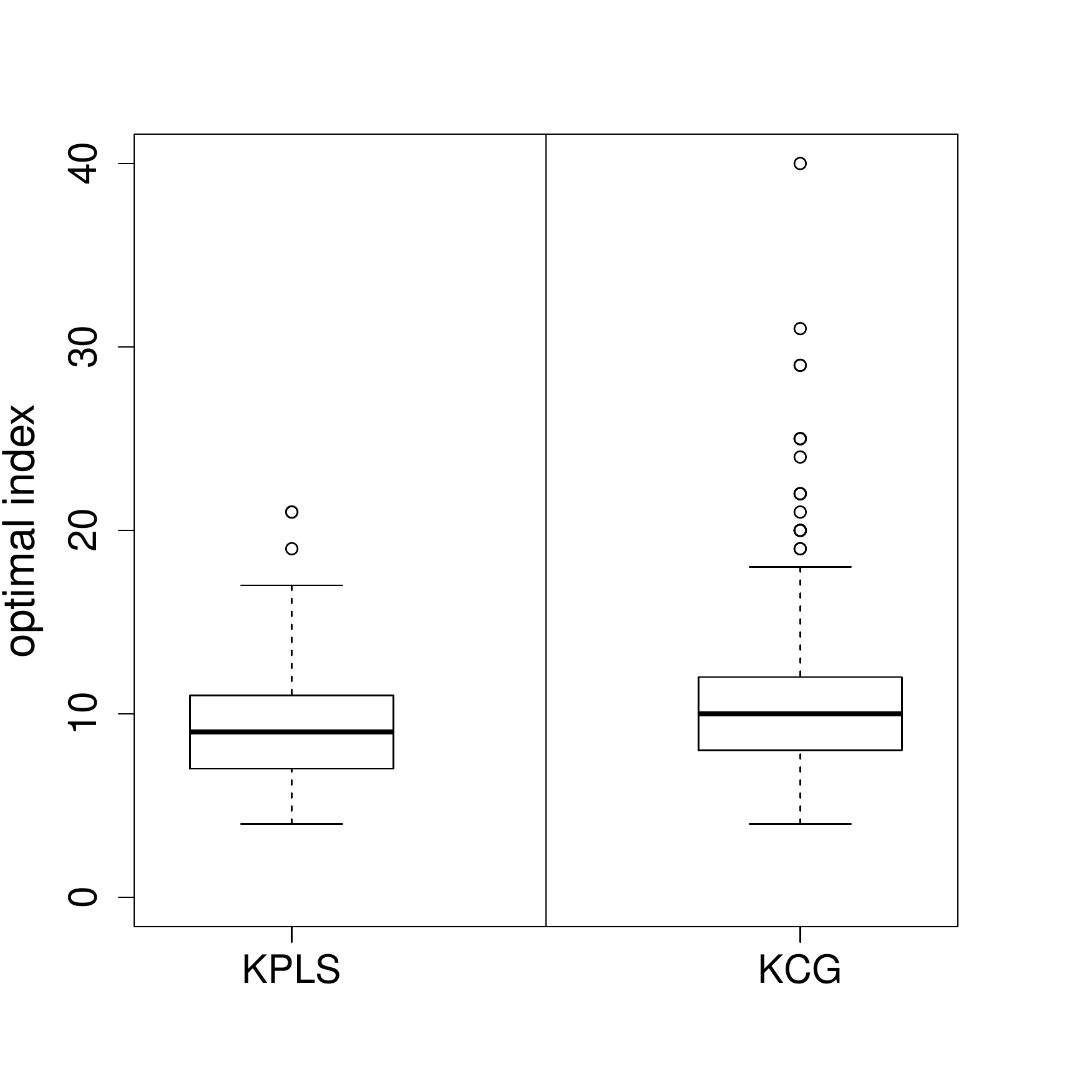}
   	\end{minipage}
   	\caption{\label{fig:a}
   		Boxplots of the optimal indices $a\in\{1,\dots,40\}$ for which the $\LTwo$-errors $\{\widehat{e}_{n,\tau}^{(j,m)}\}_{j=1}^{M}$ were attained. Kernel partial least squares is on the left of each panel and kernel conjugate gradient on the right. On the left is $n = 200$, on the right $n=1000$. The data were assumed to be independent and identically distributed.
   	}
   \end{center}
\end{figure}

We also investigated the the stopping indices $a=1,\dots,40$ for which the errors $\widehat{e}_{n,\tau}^{(j,m)}$ were attained. These are shown in Figure \ref{fig:a} for independent and identically distributed data.
It can be seen that the optimal indices for both algorithms have a rather similar behaviour. Kernel conjugate gradient stops slightly later, but overall the differences seem negligible.

\begin{figure}
   \begin{center}
   	\begin{minipage}{.45\textwidth}
	\center   	
   	\includegraphics[width=\ssize]{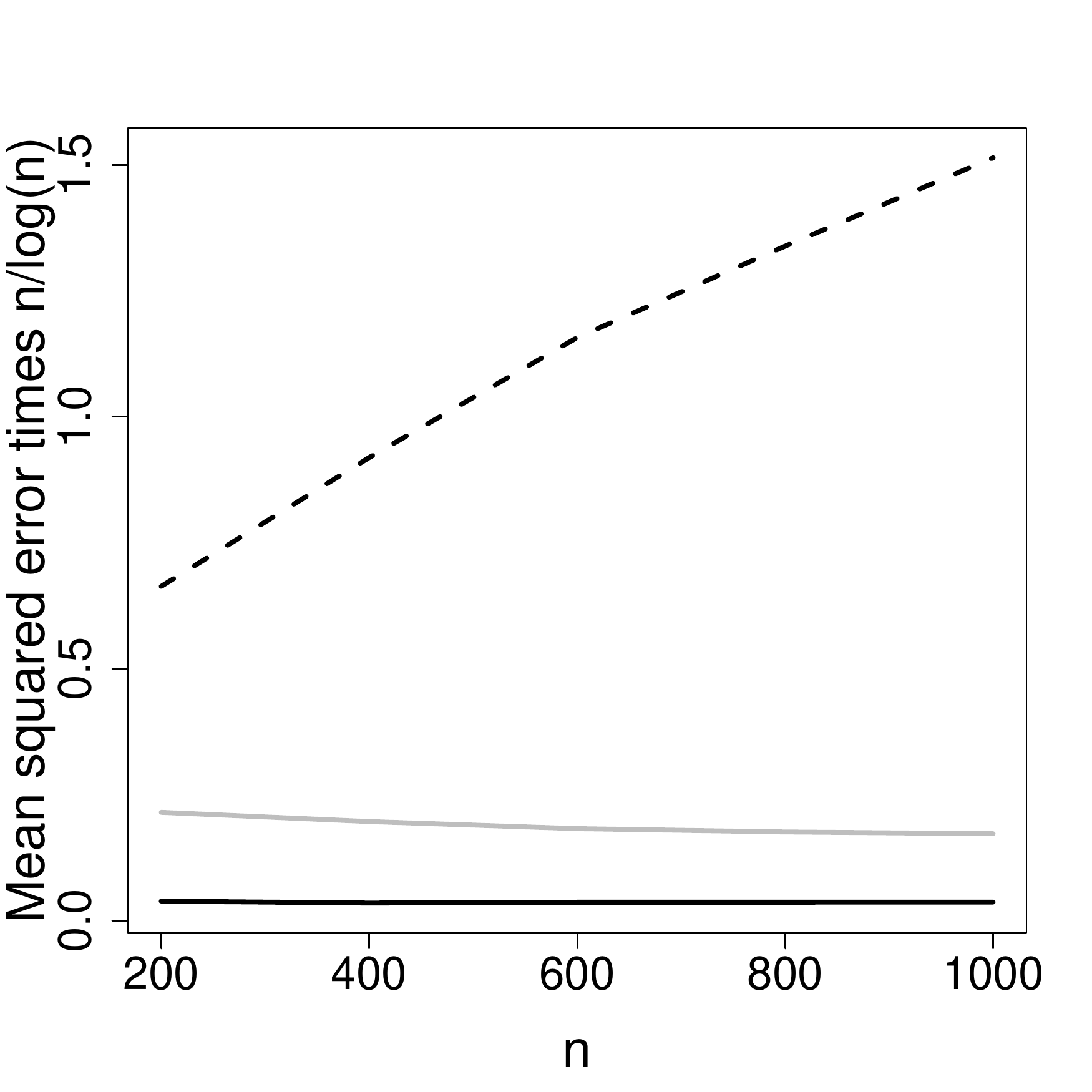}
   	\end{minipage}
   	\begin{minipage}{.45\textwidth}
	\center   	
   	\includegraphics[width=\ssize]{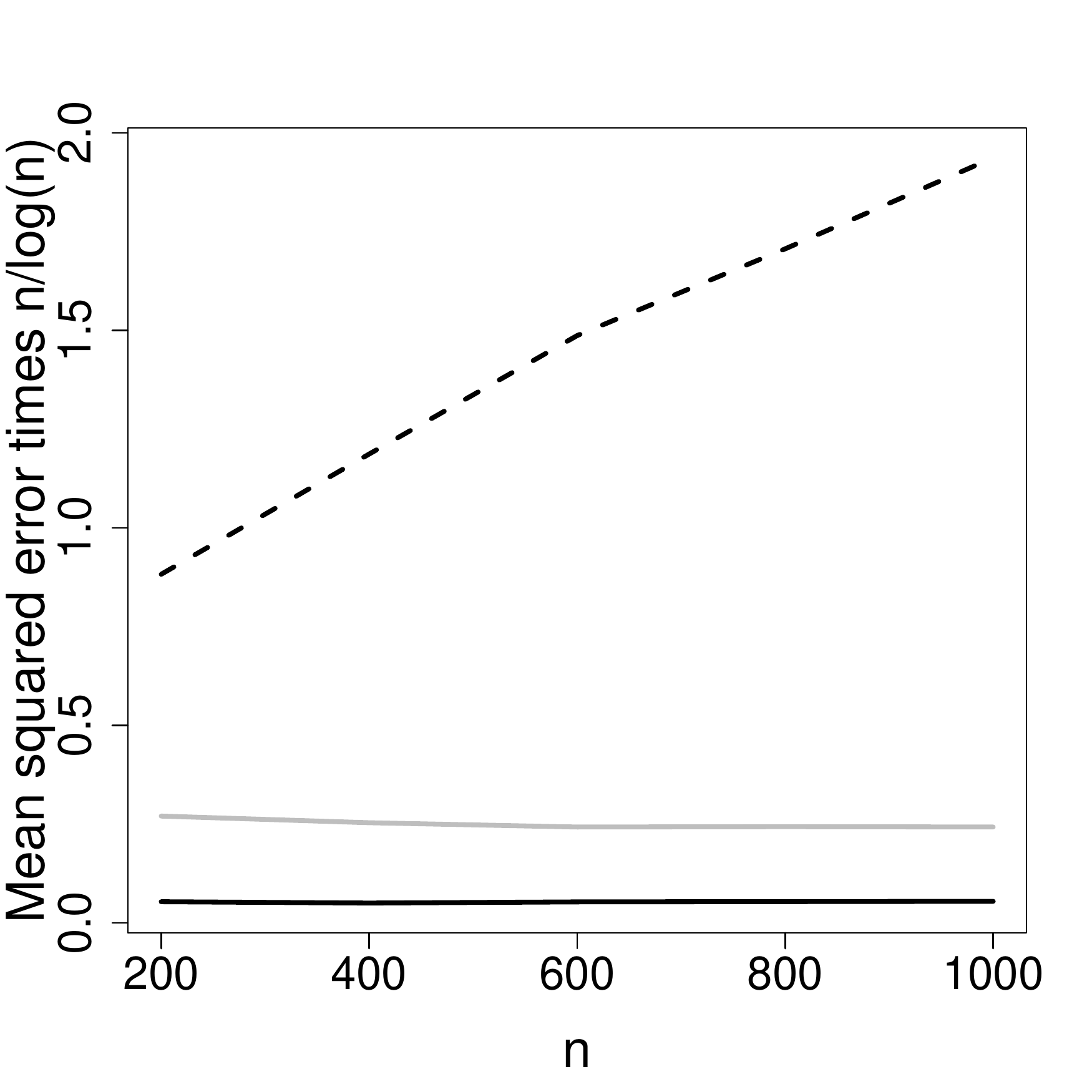}
   	\end{minipage}
   	\caption{\label{fig:lines}
   		Mean of the $\LTwo$-errors $\{\widehat{e}_{n,\tau}^{(j,m)}\}_{j=1}^{M}$ of kernel partial least squares (left) and kernel conjugate gradient (right) for $n = 200,400,\dots,1000$ multiplied by $n/\log(n)$. The solid black line is for $\tau_h = \mathbb{I}(h=0)$, the grey line for $\tau_h = 0.9^{-h}$ and the dashed black line for $\tau_h = (h+1)^{-1/4}$.
   	}
   \end{center}
\end{figure}

Figure \ref{fig:lines} shows the mean (over $j$) of the estimated $\LTwo$ errors $\{\widehat{e}_{n,\tau}^{(j,m)}\}_{j=1}^{M}$ for different $n$, $\tau$ and $m \in \{CG,PLS\}$. The errors were multiplied by $n/\log(n)$ to illustrate the convergence rates. According to Proposition \ref{prop:ed} and Corollary \ref{cor:convergence} (ii) we expect the rates for the independent and autoregressive cases to be $n^{-1}\log(n)$, which is verified by the fact that the solid black and grey lines are roughly constant.
For the long range dependent case we expect worse convergence rates which are also illustrated by the divergence of the dashed black line.

\section{Application to Molecular Dynamics Simulations}
\label{sec:protein}
The collective motions of protein atoms are responsible for its biological function and molecular dynamics simulations is a popular tool to explore this \citep{Henz07}. 

Typically, the $p \in \N$ backbone atoms of a protein are considered for the analysis with the relevant dynamics happening in time frames of nanoseconds. Although the dynamics are available exactly, the high dimensionality of the data and large number of observations can be cumbersome for regression analysis, e.g., due to the high collinearity in the columns of the covariates matrix. Many function-dynamic relationships are also non-linear \citep{Hub09}.
A further complication is the fact that the motions of different backbone atoms are highly correlated, making additive non-parametric models for the target function $f^\ast$ less suitable. 

We consider T4 Lysozyme (T4L) of the bacteriophage T4, a protein responsible for the hydrolisis of 1,4-beta-linkages in peptidoglycans and chitodextrins from bacterial cell walls.
The number of available observations is $n=4601$ and T4L consists of $p=486$ backbone atoms.

Denote with $A_{t,i} \in \R^3$ the $i$-th atom, $i=1,\dots,p$, at time $t=1,\dots,n$ and $c_i \in \R^3$ the $i$-th atom in the (apo) crystal structure of T4L. A usual representation of the protein in a regression setting is the Cartesian one, i.e., we take as the covariate $X_t = (A_{1,t}^\T,\dots,A^\T_{p,t})^\T$, $t=1,\dots,n$, see \citet{Bro83}.
The functional quantity to predict is the root mean square deviation of the protein configuration $X_t$ at time $t = 1,\dots,n$ from the (apo) crystal structure $C = (c_1^\T,\dots,c_d^\T)^\T$, i.e.,
\[
	y_t = \left\{
		p^{-1} \sum\limits_{i=1}^{p} \|X_{i,t}-C_i\|^2
	\right\}^{1/2}.
\] 
This nonlinear function was previously considered in \citet{Hub09}, where it was established that linear models are insufficient for the prediction.

\begin{figure}
   \begin{center}
	 	\begin{minipage}{.48\textwidth}
	  		\centering
	 		\includegraphics[width=\ssize]{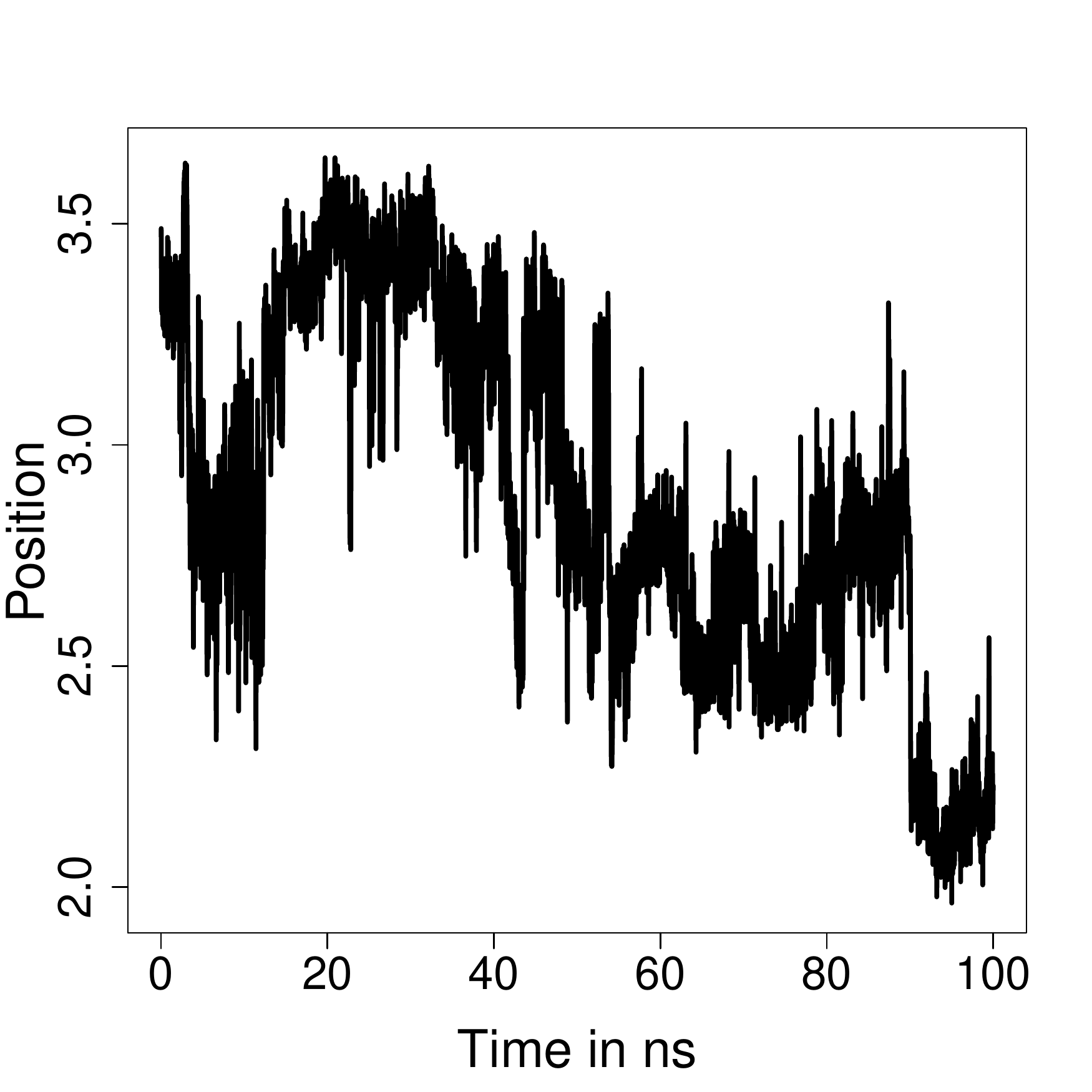}
		\end{minipage}%
		\begin{minipage}{.48\textwidth}
	  		\centering
	 		\includegraphics[width=\ssize]{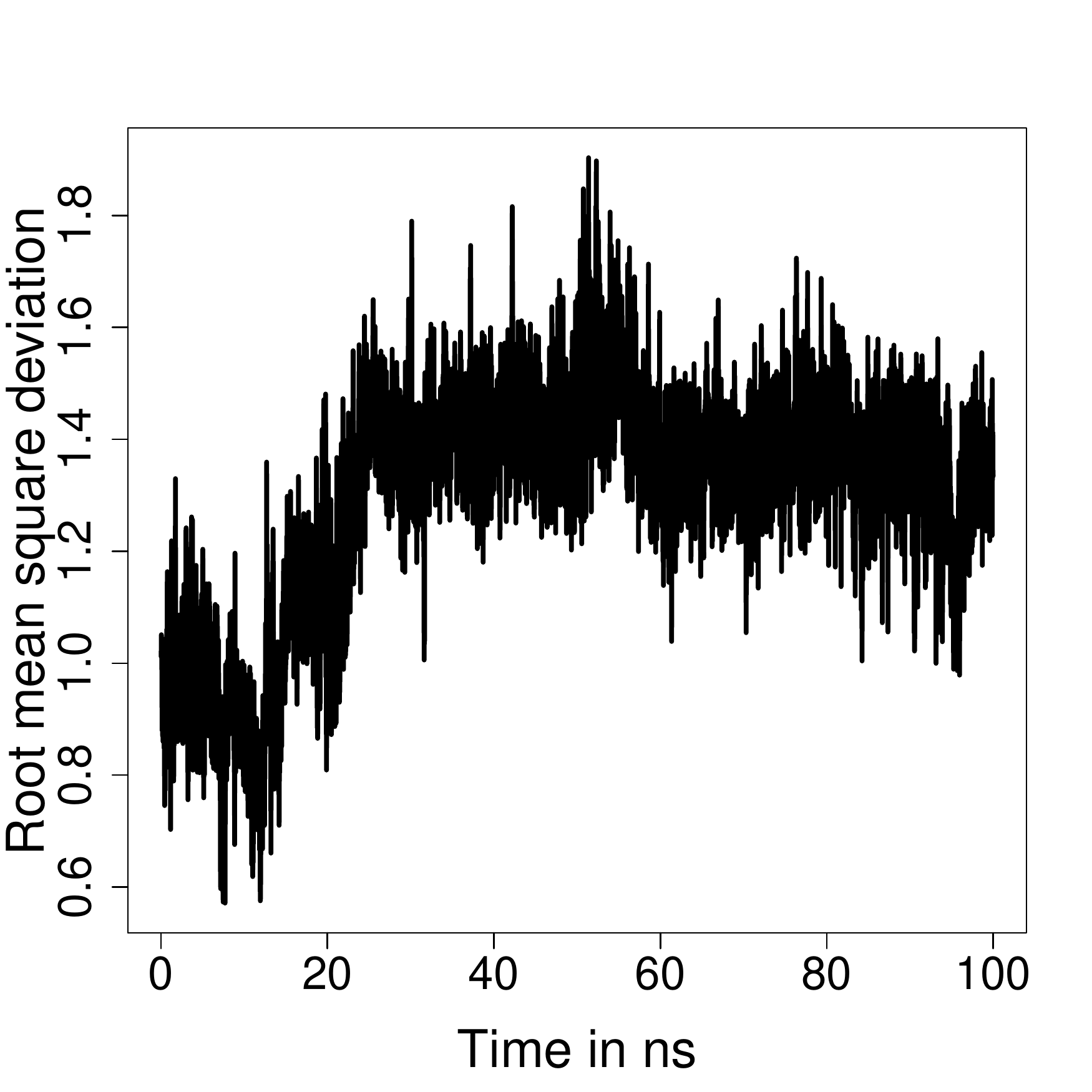}
		\end{minipage}
		\caption{\label{fig:sample.series}
			Time series of $X_{t,1}$, i.e., the first coordinate of the first atom T4L consists of (left) and the root mean squared deviation $y_t$ between the protein configuration at time $t$ and the (apo) crystal structure.
		}
	\end{center}
\end{figure}

Figure \ref{fig:sample.series} shows the time series corresponding to $X_{t,1}$ (i.e., the first coordinate of the first atom of T4L) on the left and the functional quantity $y_t$ on the right. These plots reveal certain persistent dependence over time.

Fitting autoregressive moving average models of order $(3,2)$ ($ARMA(3,2)$) to $y_t$ and $ARMA(5,2)$ to $X_{t,1}$ shows that the smallest root of their respective characteristic polynomial is close to one ($1.009$ for $y_t$ and $1.003$ for $X_{t,1}$), highlighting that we are on the border of stationarity, see, e.g., \citet{bBrock}.

\begin{figure}
   \begin{center}
	 	\begin{minipage}{.48\textwidth}
	  		\centering
	 		\includegraphics[width=\ssize]{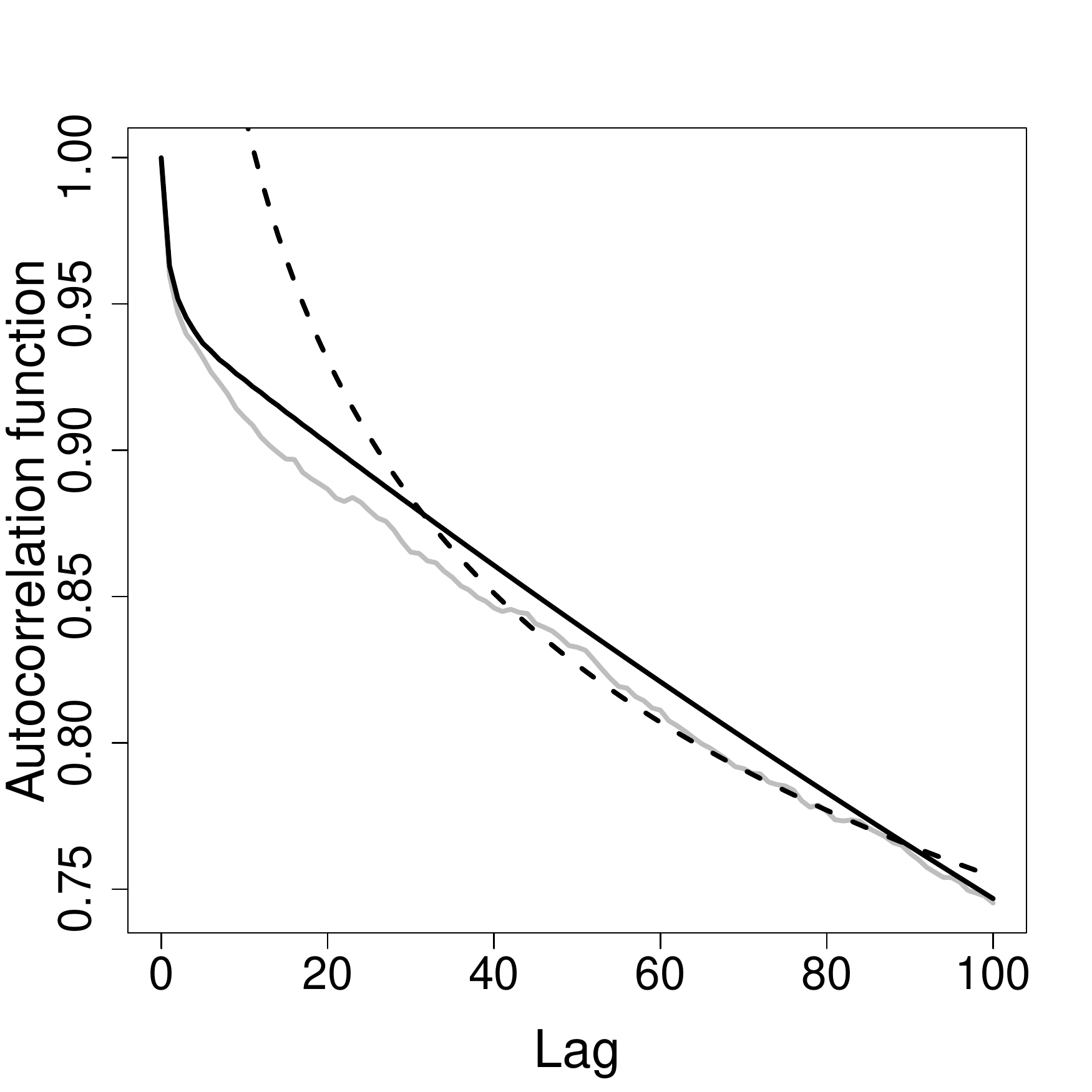}
		\end{minipage}%
		\begin{minipage}{.48\textwidth}
	  		\centering
	 		\includegraphics[width=\ssize]{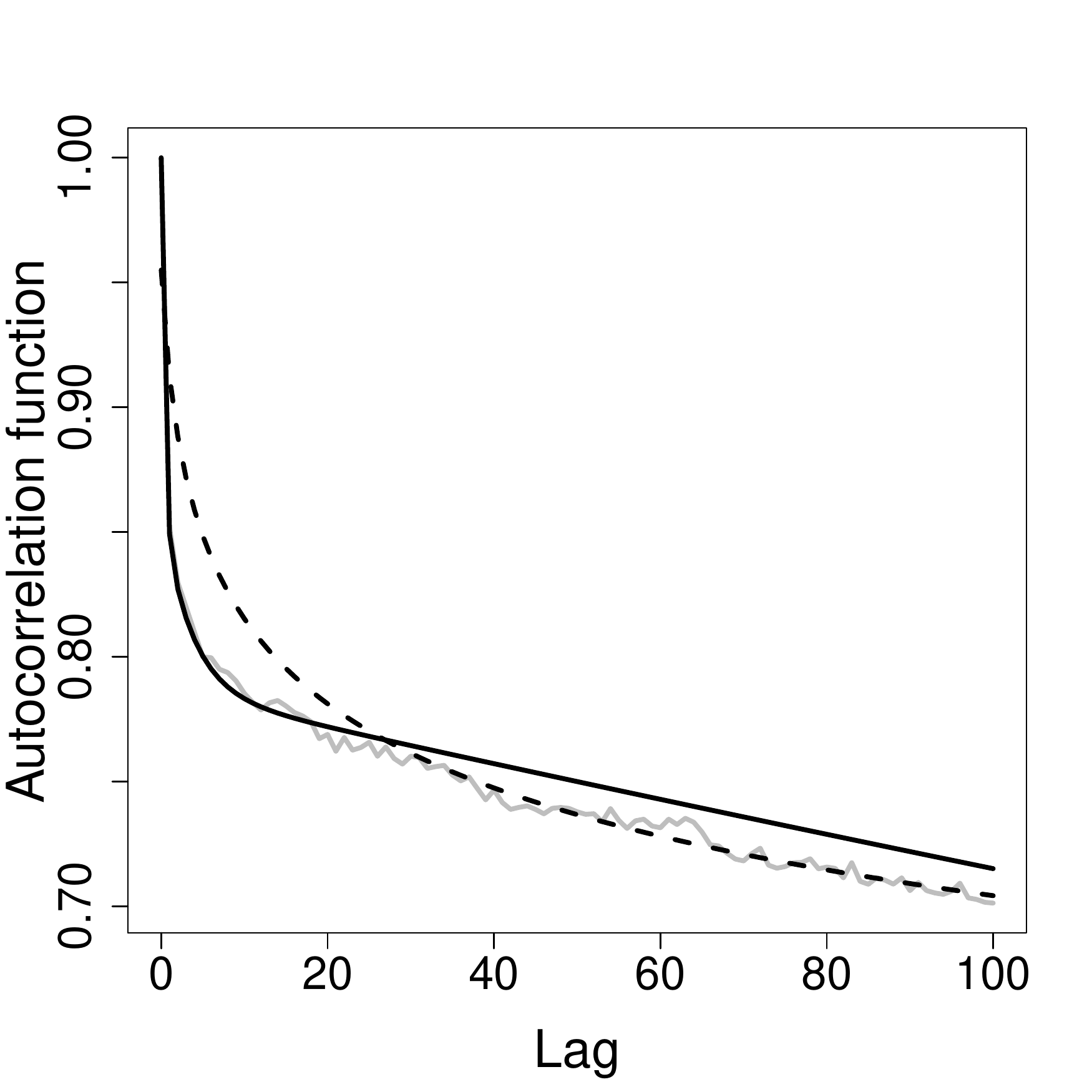}
		\end{minipage}
		\caption{\label{fig:acf}
			Autocorrelation  
			plots of $X_{t,1}$ (left) and $y_t$ (right). The estimated autocorrelation function is grey, the theoretical one of a fitted $ARMA(3,2)$ process is solid black and $\rho_h \propto (h+1)^{-q}$ for a suitable choice of $q>0$ is dashed black.
		}
	\end{center}
\end{figure}

Figure \ref{fig:acf} depicts the autocorrelation functions of $X_{t,1}$ and $y_t$, the theoretical autocorrelation function of the corresponding autoregressive moving average process and $\rho_h \propto (h+1)^{-q}$ for $q=0.134$ for $X_{t,1}$ and $q=0.066$ for $y_t$. The latter, as highlighted in Section \ref{sec:concentration}, is an autocorrelation function for a stationary long range dependent process.
These plots suggest that $X_{t,1}$ and $y_t$ follow some long-range stationary process.

We apply kernel partial least squares to this data set with the Gaussian kernel $k(x,y) = \exp(-l \|x-y\|^2)$, $x,y \in \R^{3p}$, $l>0$. The function $f$ we aim to estimate is a distance between protein configurations, so using a distance based kernel seems reasonable. Moreover, we also investigated the impact of other bounded kernels such as triangular and Epanechnikov and obtained similar results. The first $25\%$ of the data form a training set to calculate the kernel partial least squares estimator and the remaining data are used for testing.

The parameter $l>0$ is calculated via cross validation on the training set. 
In our evaluation we obtained $l = 0.0189$.

\begin{figure}[t!]
   \begin{center}
	 	\begin{minipage}{.48\textwidth}
	  		\centering
	 		\includegraphics[width=\ssize]{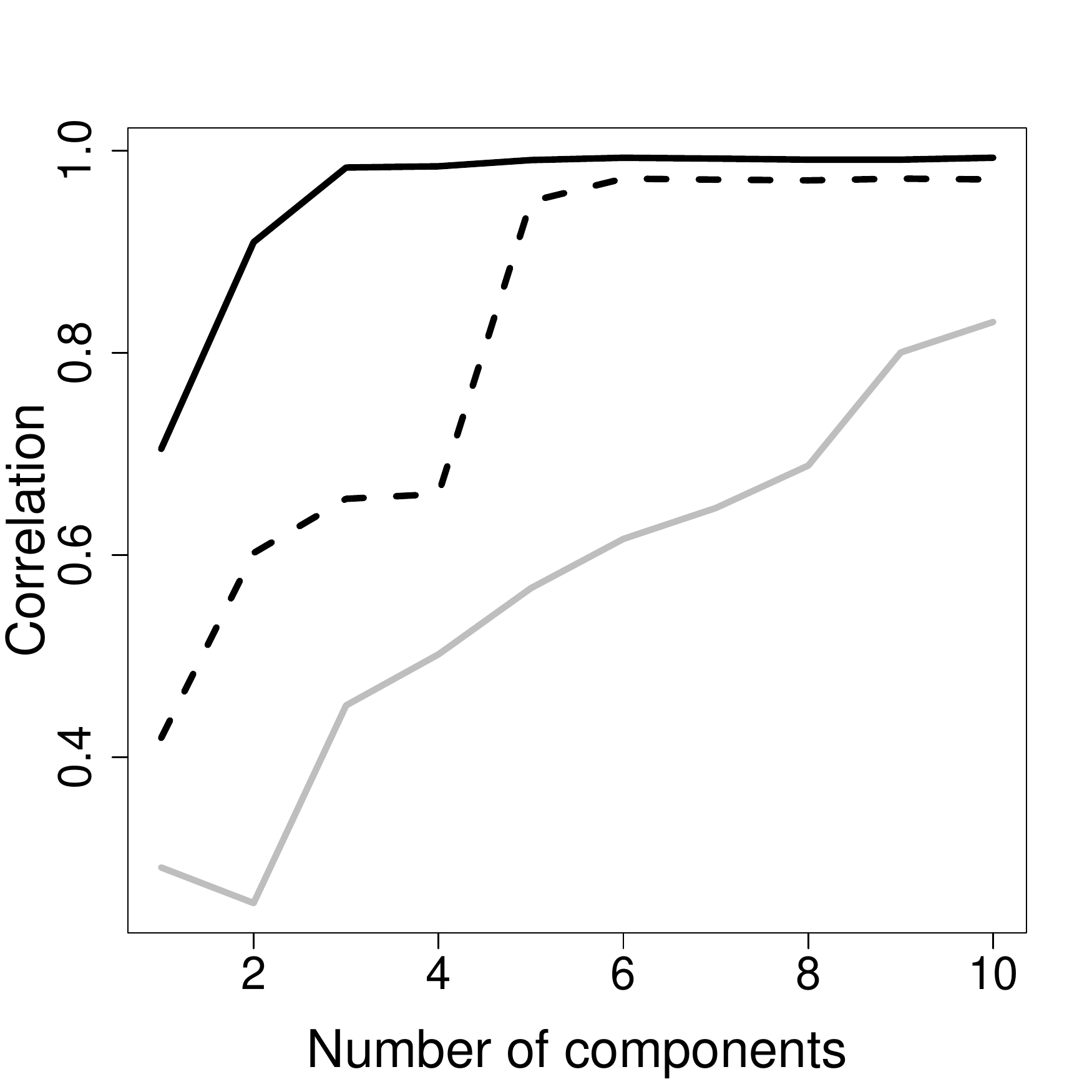}
		\end{minipage}%
		\begin{minipage}{.48\textwidth}
	  		\centering
	 		\includegraphics[width=\ssize]{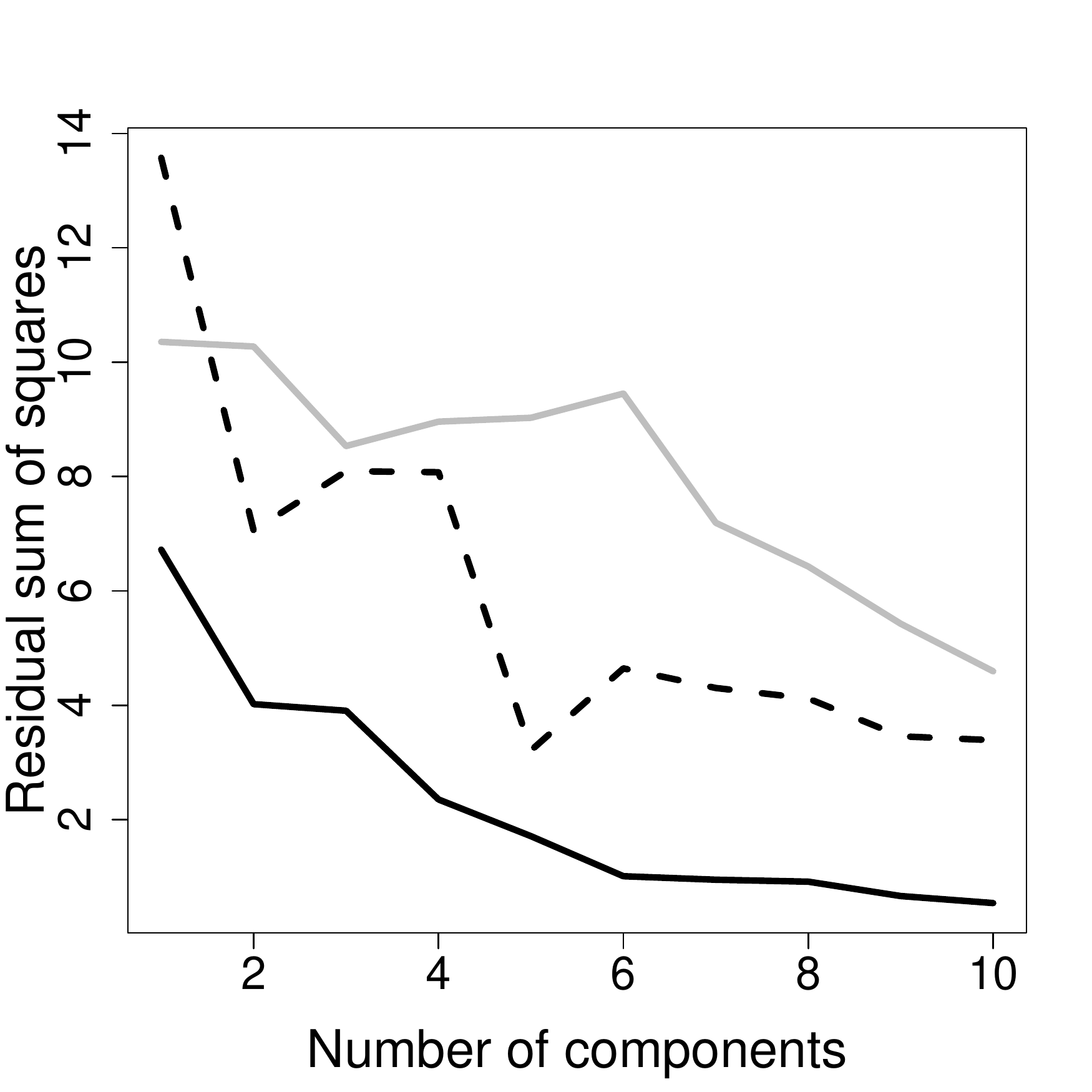}
		\end{minipage}
		\caption{\label{fig:prot}
			Correlation (left) and residual sum of squares (right) between predicted values and the observed response on the test set depending on the number of used components for kernel partial least squares (solid black), partial least squares (grey) and kernel principal component regression (dashed black).
		}
	\end{center}
\end{figure}

Figure \ref{fig:prot} compares the observed response in the test set with the prediction on the test set obtained by kernel partial least squares, kernel principal component regression and linear partial least squares.
Apparently, kernel partial least squares show the best performance and the kernel principal components algorithm is able to achieve comparable prediction with more components only. Obviously, linear partial least squares can not cope with the non-linearity of the problem. 

This application highlights that kernel partial least squares still delivers a robust prediction even when the dependence in the data is more persistent, if enough observations are available.

\bibliography{bibliography}
\biblist

\appendix
\section{Proofs}
\label{sec:proofs}
\subsection{Proof of Theorem \ref{th:kpls}}
\label{sec:int.proof1}
The proof of Theorem \ref{th:kpls} makes use of the connection between the partial least squares and the conjugate gradient algorithm. This section is structured as follows: First we will introduce the link between kernel partial least squares and kernel conjugate gradient. We will state some key facts about orthogonal polynomials and their relationship to the algorithm in Lemma \ref{lem:orth.pol}.
Then the consistency of kernel partial least squares is shown with the help of three error bounds that are obtained in Lemmas \ref{lem:error.difference} -- \ref{lem:derivative}.

With a slight abuse of notation we define $f_i = f_{\alphaest_i}$ for $i=1,\dots,n$. 
We consider the kernel partial least squares algorithm as an optimization problem
\begin{equation}
\label{ref:pls.opt}
f_i = \arg\min\limits_{g \in \mathcal{K}_i(\KEst,\bEst)} \|y - T_n g\|^2, ~~ i=1,\dots,n.
\end{equation}
This is the conjugate gradient algorithm CGNE discussed in chapter 2.2 of \citet{Hanke}.

\subsubsection{Orthogonal polynomials and some notation}
Denote with $\mathcal{P}_i$ the set of polynomials of degree at most $i = 0,\dots,n$. For functions $\psi,\phi:\R \rightarrow \R$ and $r \in \N_0$ define the inner products $[\psi,\phi]_r = \langle  \psi(\KEst)\bEst, \KEst^r \phi(\KEst) \bEst\rangle_\calH$. From the definition of the Krylov space it is immediate that every element $v \in \calK_i(\KEst,\bEst)$, $i=1,\dots,n$, can be represented by a polynomial $q \in \mathcal{P}_{i-1}$ via $v = q(\KEst) \bEst$. 

The following discussion is based on \citet{Hanke}, chapter 2. There exist two sequences of polynomials $\{p_i\}_{i=0}^n, \{q_i\}_{i=0}^n \subset \mathcal{P}_i$, such that $f_i = q_{i-1}(\KEst) \bEst$ with $q_{-1} = 0$ and $\bEst - \KEst f_i = p_i(\KEst)\bEst$. Both sequences are connected by the equation $p_i(x) = 1 - x q_{i-1}(x)$, $x \in \R$, and the polynomials $\{p_i\}_{i=0}^n$ are orthogonal with respect to $[\cdot,\cdot]_0$.

We will also consider other sequences of polynomials, namely $\{\qpol{i}{r}\}_{i=0}^n,\{\pol{i}{r}\}_{i=0}^n \subset \calP_i$, $\qpol{-1}{r} = 0$, such that $\pol{i}{r}(x) = 1 - x \qpol{i-1}{r}(x)$, $x \in \R$, and the sequence $\{\pol{i}{r}\}_{i = 0}^n$ is orthogonal with respect to $[\cdot,\cdot]_{r}$. This yields for every $r \in \N_0$ a separate conjugate gradient algorithm with solution $\fest{i}{r} = \qpol{i-1}{r}(\KEst) \bEst \in \calK_{i}(\KEst,\bEst)$ and residuals $\bEst - \KEst\fest{i}{r} = \pol{i}{r}(\KEst)\bEst$, $i=1,\dots,n$.
The $p_i^{[r]}$, $i=0,\dots,n$, $r \in \N_0$, are called residual polynomials.

As $\KEst$ is self-adjoint, positive semi-definite and the kernel is bounded by $\kappa$ we know that its spectrum is a subset of $[0,\kappa]$, see \citet{Cap07}. This also implies that $\max\{\|\KCov\|_{\calL},\|\KEst\|_\calL\} \leq \kappa$, with $\|\cdot\|_\calL$ denoting the operator norm.
The $i$ distinct roots of $\pol{i}{r}$ will be denoted by $0 < \xr{1}{i}{r} < \dots \xr{i}{i}{r} < \kappa$, $i=1,\dots,n$.

We will summarize some key facts about the orthogonal polynomials in the next lemma.
\begin{lemma}
\label{lem:orth.pol}
Let $r,s \in \N_0$ and $i=1,\dots,n$. Then we have:
\begin{itemize}
\item[(i)] 
The roots of consecutive orthogonal polynomials interlace, i.e., for $j=1,\dots,i$ it holds
\[
	0 < \xr{j}{i+1}{r} < \xr{j}{i}{r}<\xr{j}{i}{r+1}<\xr{j+1}{i+1}{r}<\xr{j+1}{i}{r} < \dots <\xr{i}{i}{r+1} < \xr{i+1}{i+1}{r} < \kappa,
\]
\item[(ii)]
the optimality property 
$[\pol{i}{1},\pol{i}{1}]_{0}^{1/2} = \|\bEst - \KEst \fest{i}{1}\|_\calH \leq \|\bEst - \KEst h\|_\calH$ holds for all $h \in \calK_{i}(\KEst,\bEst)$,
\item[(iii)] on $x \in [0,\xr{1}{i}{r}]$ it holds
$0 \leq \pol{i}{r}(x) \leq 1$ and $\qpol{i}{r}(x) \leq 
\left|\left(\pol{i}{r}\right)'(0)\right|
$,
\item[(iv)]
$
\pol{n}{r} = \pol{n}{s},
$
\item[(v)] 
$\left(\pol{i}{r}\right)'(0) = - \sum_{j=1}^i \left(\xr{j}{i}{r}\right)^{-1}$,
\item[(vi)]
for $r \geq 1$ define $\phi_i(x) = \pol{i}{r}(x)\left(\xr{1}{i}{r}\right)^{1/2} \left( \xr{1}{i}{r}- x\right)^{-1/2}$, $x \in [0,\xr{1}{i}{r}]$, $i=1,\dots,n$. Then it holds for $u \geq 0$ that $x^u \phi^2_i(x) \leq u^u \left|\left(\pol{i}{r}\right)'(0)\right|^{-u}$ with the convention $0^0 = 1$.
\end{itemize}
\end{lemma}
{\textit{Proof}}: (i) See \citet{Hanke}, Corollary 2.7.

(ii) See \citet{Hanke}, Proposition 2.1.

(iii) Due to part (i) we know that all $i$ roots of the polynomial $\pol{i}{r}$ are contained in $(0,\kappa)$. Furthermore $\pol{i}{r}(0) = 1 - 0\qpol{i}{r} = 1$. Thus $\pol{i}{r}$ is convex and falling in $[0,\xr{1}{i}{r}]$ and the first assertion follows.

Because of the convexity of $\pol{i}{r}$ on $[0,\xr{1}{i}{r}]$ we get $\qpol{i}{r} (x) = x^{-1}\{1-\pol{i}{r}(x)\} \leq \left|\left(\pol{i}{r}\right)'(0)\right|$.

(iv) See the discussion in \citet{Hanke} preceding Proposition 2.1 and use the facts that
$\bEst \in \mathrm{range}(\KEst)$ and $\KEst$ is an operator of rank $n$.

(v) Write $\pol{i}{r}(x) = \prod_{j=1}^i(1 - x/ \xr{j}{i}{r})$, $x \in [0,\kappa]$, and the result is immediate.

(vi) See equation (3.10) in \citet{Hanke}.
\eop


We denote for $x\geq 0$ by $\Proj_x$ the orthogonal projection operator on the eigenspace corresponding to the eigenvalues of $\KEst$ that are smaller or equal $x$ and $\Proj_x^\perp = I_\calH - P_x$ with $I_\calH:\calH \rightarrow \calH$ being the identity operator. 


\subsubsection{Preparation for the proof}
An important technical result that will be useful in the upcoming proof is
\begin{lemma}
\label{lem:op.inequality}
Let $B,C:\calH \rightarrow \calH$ be two positive semi-definite, self-adjoint operators with $\max\{\|B\|_\calL,\|C\|_\calL\} \leq \kappa$. Then it holds for any $r \geq 0$ with $\zeta = \max\{r-1,0\}$
\[
	\|B^r - C^r\|_\calL \leq (\zeta+1) \kappa^\zeta \|B - C\|_\calL^{r-\zeta}.
\]
\end{lemma}
{\textit{Proof}}: See \citet{Blan10b}, Lemma A.6. \eop

For the remainder of the proof we assume that we are on the set where it holds with probability at least $1-\nu$, $\nu \in (0,1]$, that $\|\KEst - \Sgs\|_\mathcal{L} \leq C_\delta(\nu) \gamma_n$ and $\|\bEst - \Sgs f\| \leq C_\epsilon(\nu) \gamma_n$ for a sequence $\{\gamma_n\}_n$ converging to zero and constants $C_\delta=C_\delta(\nu),C_\epsilon=C_\epsilon(\nu)>0$.

With Lemma 2.4 in \citet{Hanke} we see that the stopping iteration (\ref{eq:stopping}) can also be expressed as 
\begin{equation}
\label{eq:alt.stop}
	\aopt = \min\left\{1 \leq a \leq n: \|\KEst \fest{a}{1}- \bEst\|_\calH \leq \Ctau \gamma_n \right\},
\end{equation}
i.e., we stop the kernel partial least squares algorithm when a discrepancy principle for $\fest{a}{1}$ holds.

It is easy to see that from \ref{eq:source} it follows for $r \geq 1/2$ that 
\begin{enumerate}[label={(SH})]
\item
\label{eq:source.H}
There exist $\mu \geq 0$, $R >0$ and $u \in \calH$ such that $f = \KCov^{\mu} u$ and $\|u\|_\calH \leq R$.
\end{enumerate}
This condition is known as the H\"older source condition with $\mu = r-1/2$. 

Recall that $\calH \subseteq \LTwo$ and $T:\calH \rightarrow \LTwo$ is the change of space operator.
Using the fact that $T$, $T^\ast$ are adjoint operators, $f_\aopt = T f_\aopt$ and $f^\ast = T f$ for $r \geq 1/2$ we see
\begin{align*}
	\|f_\aopt - f^\ast\|_2 = \|T(f_\aopt - f)\|_2 = \langle S(f_\aopt - f), f_\aopt - f\rangle_\calH = \|S^{1/2}(f_\aopt - f)\|_\calH.
\end{align*}
An application of Lemma \ref{lem:op.inequality} yields
\begin{align}
	&\|f_\aopt - f^\ast\|_2 =\|\Sg (f_\aopt - f)\|_\calH \leq \|\Sg (f_\aopt - \fest{\aopt}{1})\|_\calH + \|\Sg (\fest{\aopt}{1} - f)\|_\calH \notag \\
	& \leq 
	C_\delta^{1/2}\gamma_n^{1/2} \left( \|f_\aopt - \fest{\aopt}{1}\|_\calH + \|\fest{\aopt}{1} - f\|_\calH
	\right) + \|\KEst^{1/2}(f_\aopt - \fest{\aopt}{1})\|_\calH  
	+ \|\KEst^{1/2}(\fest{\aopt}{1}-f)\|_\calH.  \label{eq:main.inequality}
\end{align} 

The following lemmas will deal with bounding the quantities in (\ref{eq:main.inequality}).

\begin{lemma}
\label{lem:error.difference}
Assume $C_x \in (0,1]$ such that $x_\ast =  (C_x\gamma_n)^{1/(\mu+1)} < \xr{1}{\aopt-1}{1}$ and $\Ctau > C_\epsilon + C_x R + C_\delta (\mu+1) \kappa^\mu R$.
Under the conditions of the theorem it holds $\mu \geq 0$
\begin{align*}
	\|f_\aopt - \fest{\aopt}{1}\|_\calH & \leq 
	\gamma^{\mu/(\mu+1)}_n
	\frac{\Ctau}
	{C_x^{1/(\mu+1)}
	\left[
	1-	\Ctau^{-1}
	\{C_\epsilon + C_x R + C_\delta (\mu+1) \kappa^\mu R\}
	\right]^2}\\
	\|\KEst^{1/2}(f_\aopt - \fest{\aopt}{1})\|_\calH & \leq 
	\gamma^{(2\mu+1)/(2\mu+2)}_n
	\frac{\Ctau}
	{C_x^{1/(2\mu+2)}
	\left[
	1-	\Ctau^{-1}
	\{C_\epsilon + C_x R + C_\delta (\mu+1) \kappa^\mu R\}
	\right]}.
\end{align*}
\end{lemma}
{\it{Proof:}}
If the inner products $[\cdot,\cdot]_0$ and $[\cdot,\cdot]_1$ are the same the proof is done because both polynomial sequences are identical.  

We now observe that we have for $\aopt = n$ due to Lemma \ref{lem:orth.pol} (iv) $q_{n-1}(x) - \qpol{n-1}{1}(x) = x^{-1}\{\pol{n}{1}(x)-p_{n}(x)\} = 0$, i.e., $\|f_{\aopt}-\fest{\aopt}{1}\|_\calH = 0$ and $\|\KEst^{1/2}(f_{\aopt}-\fest{\aopt}{1})\|_\calH = 0$ and the proof is done.

If the inner products differ and we have $0<\aopt<n$ it holds $f_\aopt \neq \fest{\aopt}{1}$.

Proposition 2.8 in \citet{Hanke} can now be applied for $0 < \aopt < n$ and yields  $q_{\aopt-1}(x) - \qpol{\aopt-1}{1}(x) = x^{-1}\{\pol{\aopt}{1}(x)-p_{\aopt}(x)\} = \theta_{\aopt} \pol{\aopt-1}{2}(x)$, $x\geq 0$, with $\theta_{\aopt} = (\pol{\aopt}{1})'(0)-(\pol{\aopt}{0})'(0)>0$. We get
$f_\aopt - \fest{\aopt}{1} = q_{\aopt-1}(\KEst)\bEst - \qpol{\aopt-1}{1}(\KEst) \bEst = \theta_{\aopt} \pol{\aopt-1}{2}(\KEst) \bEst$. 

Proposition 2.9 in \citet{Hanke} yields
$\theta_{\aopt} = \polnorm{\pol{\aopt-1}{2}}{1}^{-1} \polnorm{\pol{\aopt}{1}}{0}$. The optimality property of $\fest{\aopt}{1}$ in Lemma \ref{lem:orth.pol} (ii) shows that 
\begin{equation}
\label{eq:pol.opt}
\|\bEst - \KEst \fest{\aopt}{1}\|_\calH = \|\pol{\aopt}{1}(\KEst)\bEst\|_\calH = \polnorm{\pol{\aopt}{1}}{0}^{1/2} \leq \polnorm{\pol{\aopt-1}{2}}{0}^{1/2}.
\end{equation}

Combining these results yields
\begin{equation}
\label{eq:est.difference}
 	\|f_\aopt - \fest{\aopt}{1}\|_\calH = \frac{\polnorm{\pol{\aopt}{1}}{0}}{\polnorm{\pol{\aopt-1}{2}}{1}} \polnorm{\pol{\aopt-1}{2}}{0}^{1/2} 
 	\leq 
	\frac{\polnorm{\pol{\aopt-1}{2}}{0}}{\polnorm{\pol{\aopt-1}{2}}{1}}\|\pol{\aopt}{1}(\KEst) \bEst\|_\calH.
\end{equation}

Recall that $\xr{1}{\aopt-1}{2}$ denotes the first root of $\pol{\aopt-1}{2}$.
It holds for any $0 \leq x \leq \xr{1}{\aopt-1}{2}$ that $0 \leq \pol{\aopt-1}{2}(x) \leq 1$, see Lemma \ref{lem:orth.pol} (iii), and thus
\begin{align*}
	\polnorm{\pol{\aopt-1}{2}}{0}^{1/2} &\leq
	\|\Proj_x \pol{\aopt-1}{2}(\KEst) \{\bEst - S f + Sf\}\|_\calH + \|\Proj_x^\perp \pol{\aopt-1}{2}(\KEst)\bEst\|_\calH\\
	&\leq C_\epsilon\gamma_n + \|\Proj_x \pol{\aopt-1}{2}(\KEst) \KCov^{\mu+1} u\|_\calH + x^{-1/2} \|\Proj_x^\perp \KEst^{1/2} \pol{\aopt-1}{2}(\KEst) \bEst\|_\calH\\
	& \leq C_\epsilon\gamma_n + x^{\mu+1} R  + \|\Proj_x \pol{\aopt-1}{2}(\KCov^{\mu+1}-\KEst^{\mu+1})u \|_\calH + x^{-1/2} \polnorm{\pol{\aopt-1}{2}}{1}^{1/2}.
\end{align*}
In the second inequality \ref{eq:source.H} with $\mu \geq 0$ was applied.

By assumption $x_\ast =  (C_x\gamma)^{1/(\mu+1)} \leq \xr{1}{\aopt-1}{1} < \xr{1}{\aopt-1}{2}$ due to the interlacing property of the roots of the polynomials $\pol{i}{r}$, $i=1,\dots,n$, $r \in \N_0$, see Lemma \ref{lem:orth.pol} (i).

Using Lemma \ref{lem:op.inequality} we get $\|\KCov^{\mu+1}-\KEst^{\mu+1}\|_\calL \leq (\mu+1) \kappa^\mu C_\delta\gamma_n$ and setting $x=x_\ast$ we get
\begin{align}
 \polnorm{\pol{\aopt-1}{2}}{0}^{1/2} &\leq
 C_\epsilon \gamma_n + x_\ast^{\mu+1} R  + C_\delta \gamma_n(\mu+1) \kappa^\mu R + x_\ast^{-1/2} \polnorm{\pol{\aopt-1}{2}}{1}^{1/2} \notag\\
 & = \gamma_n \left\{
 	C_\epsilon + C_x R + C_\delta (\mu+1) \kappa^\mu R
 \right\} + x_\ast^{-1/2} \polnorm{\pol{\aopt-1}{2}}{1}^{1/2}.\label{eq:gamma.bound}
\end{align}
Due to (\ref{eq:alt.stop}) and (\ref{eq:pol.opt}) we have additionally
$
	\Ctau \gamma_n \leq  \|\KEst \fest{\aopt-1}{1}- \bEst\|_\calH = \|\pol{\aopt-1}{1}(\KEst) \bEst\|_\calH \leq \polnorm{\pol{\aopt-1}{2}}{0}^{1/2}.
$

Plugging this into (\ref{eq:gamma.bound}) yields
\[
	\polnorm{\pol{\aopt-1}{2}}{0}^{1/2} \leq 
	\Ctau^{-1}
	\{C_\epsilon + C_x R + C_\delta (\mu+1) \kappa^\mu R\}
	\polnorm{\pol{\aopt-1}{2}}{0}^{1/2}
	+ 
	x_\ast^{-1/2} \polnorm{\pol{\aopt-1}{2}}{1}^{1/2},
\]
or equivalently with $x_\ast =  (C_x\gamma_n)^{1/(\mu+1)}$
\begin{equation}
\label{eq:twosidepol}
	\polnorm{\pol{\aopt-1}{2}}{0}^{1/2} \leq 
	\gamma^{-1/(2\mu+2)}_n	
	\frac{\polnorm{\pol{\aopt-1}{2}}{1}^{1/2}}
	{C_x^{1/(2\mu+2)}
	\left[
	1-	\Ctau^{-1}
	\{C_\epsilon + C_x R + C_\delta (\mu+1) \kappa^\mu R\}
	\right]  }
	 ,
\end{equation}
where by assumption $\Ctau > C_\epsilon + C_x R + C_\delta (\mu+1) \kappa^\mu R$ and $x_\ast = (C_x\gamma)^{1/(\mu+1)}$.

Combining (\ref{eq:est.difference}), (\ref{eq:twosidepol}) and $\|\pol{\aopt}{1}(\KEst) \bEst\|_\calH \leq \Ctau \gamma_n$ due to the stopping index (\ref{eq:alt.stop}) yields
\begin{align*}
	\|f_\aopt - \fest{\aopt}{1}\|_\calH &\leq 
	\gamma^{-1/(\mu+1)}_n
	\frac{
	\|\pol{\aopt}{1}(\KEst) \bEst\|_\calH
	}
	{C_x^{1/(\mu+1)}
	\left[
	1-	\Ctau^{-1}
	\{C_\epsilon + C_x R + C_\delta (\mu+1) \kappa^\mu R\}
	\right]^2  }\\
	& \leq 
	\gamma^{\mu/(\mu+1)}_n
	\frac{\Ctau}
	{C_x^{1/(\mu+1)}
	\left[
	1-	\Ctau^{-1}
	\{C_\epsilon + C_x R + C_\delta (\mu+1) \kappa^\mu R\}
	\right]^2}.
\end{align*}
For the second part of the proof we derive in the same way as (\ref{eq:est.difference})
\[
	\|\KEst^{1/2}(f_\aopt - \fest{\aopt}{1})\|_\calH \leq 
	\frac{\polnorm{\pol{\aopt-1}{2}}{0}^{1/2}}
	{\polnorm{\pol{\aopt-1}{2}}{1}^{1/2}}
	\|\pol{\aopt}{1}(\KEst)\bEst\|_\calH.
\]
Using (\ref{eq:twosidepol}) and $\|\pol{\aopt}{1}(\KEst) \bEst\|_\calH \leq \Ctau \gamma_n$ gives
\[
	\|\KEst^{1/2}(f_\aopt - \fest{\aopt}{1})\|_\calH \leq 
	\gamma^{(2\mu+1)/(2\mu+2)}_n
	\frac{\Ctau}
	{C_x^{1/(2\mu+2)}
	\left[
	1-	\Ctau^{-1}
	\{C_\epsilon + C_x R + C_\delta (\mu+1) \kappa^\mu R\}
	\right]},
\]
finishing the proof. \eop

\begin{lemma}
\label{lem:secondlemma}
For any $i =1,\dots,n$ and any $0 < x \leq \xr{1}{i}{1}$ we have under the conditions of the theorem for $\mu \geq 1$
\begin{align*}
	\|f - \fest{i}{1}\|_\calH &\leq 
		R \left\{
		 	x^\mu  + C_\delta  \mu \kappa^{\mu-1}\gamma_n
		 \right\}	
	+ 
	x^{-1} \left\{
	\|\KEst \fest{i}{1} - \bEst\|_\calH + (C_\epsilon + C_\delta \kappa^\mu R)\gamma_n
	\right\}\\
	&+
	(C_\epsilon + C_\delta \kappa^\mu R)\gamma_n | (\pol{i}{1})'(0)|,
	\\
	\|\KEst^{1/2}(f - \fest{i}{1})\|_\calH &\leq 
	R \left\{
		 	x^{\mu+1/2}  + x^{1/2}C_\delta \mu \kappa^{\mu-1}\gamma_n
		 \right\}	\\
	&+ 
	x^{-1/2} \left\{
	\|\KEst \fest{i}{1} - \bEst\|_\calH + (C_\epsilon + C_\delta \kappa^\mu R)\gamma_n
	\right\}\\
	&+
	x^{1/2}(C_\epsilon + C_\delta \kappa^\mu R)\gamma_n | (\pol{i}{1})'(0)|.
\end{align*}
\end{lemma}
{\it{Proof}:} Denote $\bar{f}_i = \qpol{i-1}{1}(\KEst) \KEst f$ and consider
\begin{equation}
\label{eq:eqsplit}
\|f - \fest{i}{1}\|_\calH \leq 
	\|\Proj_x(f - \bar{f}_i)\|_\calH + \|\Proj_x(\bar{f}_i - \fest{i}{1})\|_\calH + 
	\|\Proj_x^\perp ( f- \fest{i}{1})\|_\calH.
\end{equation}
The first term of (\ref{eq:eqsplit}) can be bound by an application of Lemma \ref{lem:op.inequality} and \ref{eq:source.H} with $\mu \geq 1$
\begin{align*}
		\|\Proj_x(f - \bar{f}_i)\|_\calH &=
		\|\Proj_x\{I - \qpol{i-1}{1}(\KEst)\KEst\}f \|_\calH 
		=
		 \|\Proj_x\pol{i}{1}(\KEst) f \|_\calH 
		 = 
		 \|\Proj_x\pol{i}{1}(\KEst) \KCov^\mu u \|_\calH  \\
		 &\leq 		 	
		 		\|\Proj_x\pol{i}{1}(\KEst) \KEst^\mu u \|_\calH
		 		+ \|\Proj_x\pol{i}{1}(\KEst) (\KCov^\mu - \KEst^\mu) u \|_\calH
		 	\\
		 & \leq 
		 R \left\{
		 	x^\mu + C_\delta \mu \kappa^{\mu-1}\gamma_n
		 \right\}.
\end{align*}
In the last inequality we used that on $0 \leq x \leq \xr{1}{i}{1}$ we have $0 \leq \pol{i}{1}(x) \leq 1$.

For the second term of (\ref{eq:eqsplit}) we use Lemma \ref{lem:orth.pol} (iii) $\qpol{i}{1}(x) \leq |(\pol{i}{1})'(0)|$ on $x \in [0,\xr{1}{i}{1}]$. This yields
\begin{align*}
	\|\Proj_x(\fest{i}{1} - \bar{f}_i)\|_\calH 
	& = 
	\|\Proj_x \qpol{i}{1}(\KEst) (\KEst f - \bEst)\|_\calH\\
	& \leq 
	\|\Proj_x \qpol{i}{1}(\KEst) (\KCov f - \bEst)\|_\calH + \|\Proj_x \qpol{i}{1}(\KEst)(\KEst - \KCov) f\|_\calH\\
	& \leq 
	(C_\epsilon + C_\delta \kappa^\mu R)\gamma_n \pdiff{i}.
\end{align*}
Finally, we have
\begin{align*}
\|\Proj_x^\perp ( f- \fest{i}{1})\|_\calH 
&\leq 
x^{-1} \|\Proj_x^\perp \KEst( f- \fest{i}{1})\|_\calH\
 \leq x^{-1}
\left\{
	\|\KEst \fest{i}{1} - \bEst\|_\calH + 
	\|\Proj_x(\bEst - \KEst f)\|_\calH
\right\}
\\
&\leq 
x^{-1} \left\{
	\|\KEst \fest{i}{1} - \bEst\|_\calH + (C_\epsilon + C_\delta \kappa^\mu R) \gamma_n
\right\}
\end{align*}
and thus the first inequality is proven.

For the second inequality we use
\begin{equation*}
\|\KEst^{1/2}(f - \fest{i}{1})\|_\calH \leq 
	\|\Proj_x \KEst^{1/2}(f - \bar{f}_i)\|_\calH + \|\Proj_x \KEst^{1/2}(\bar{f}_i - \fest{i}{1})\|_\calH + 
	\|\Proj_x^\perp \KEst^{1/2}( f- \fest{i}{1})\|_\calH.
\end{equation*}
In the same way as before we derive bounds for the three terms:
\begin{align*}
\|\Proj_x \KEst^{1/2}(f - \bar{f}_i)\|_\calH & \leq 
x^{1/2} C_\delta \mu \kappa^{\mu-1} R\gamma_n + x^{\mu+1/2} R,\\
\|\Proj_x \KEst^{1/2}(\bar{f}_i - \fest{i}{1})\|_\calH &\leq
x^{1/2} ( C_\epsilon + C_\delta R \kappa^\mu)\gamma_n \pdiff{i},\\
\|\Proj_x^\perp \KEst^{1/2}( f- \fest{i}{1})\|_\calH  &\leq 
x^{-1/2} \{ \|\KEst \fest{i}{1} - \bEst\|_\calH + (C_\epsilon + C_\delta \kappa^\mu R)\gamma_n\},
\end{align*}
completing the proof. \eop

\begin{lemma}
\label{lem:derivative}
Assume that $C_x \in (0,1]$ is such that $x_\ast =  (C_x\gamma)^{1/(\mu+1)} < \xr{1}{\aopt-1}{1}$ and $\Ctau > C_\epsilon + C_x R + C_\delta (\mu+1) \kappa^\mu R$.
Under the conditions of the theorem it holds for $\mu \geq 0$
\begin{align*}
	\pdiff{\aopt}  &\leq 
	\gamma^{-1/(\mu+1)}_n \left[
	C_x^{-1/(\mu+1)}	
	\left\{
		1  - \frac{C_\epsilon + C_x R + C_\delta (\mu+1)\kappa^\mu R}{\Ctau}
	\right\}^{-2}
	\right.\\
	&+ \left.
	\left\{
	  	\frac{(2\mu+2)^{\mu+1} R}
	{
		\Ctau - 	C_\delta (\mu+1) \kappa^\mu R + C_\epsilon
	}
	  \right\}^{1/(\mu+1)}
	\right]
\end{align*}
\end{lemma}
{\it{Proof:}}
The proof is done in two steps by using the inequality $\pdiff{\aopt} \leq \pdiff{\aopt-1} + \left|
	\left(
		\pol{\aopt}{1}
	\right)'(0)
	-
	\left(
		\pol{\aopt-1}{1}
	\right)'(0)
\right|$.

Consider first $\aopt >1$. 

We will bound $\|\KEst \fest{\aopt-1}{1} - \bEst\|_\calH$ from above. Define $z =\xr{1}{\aopt-1}{1}$ and $\phi_{i}(x) = \pol{i}{1}(x) (z-x)^{-1/2} z^{1/2}$, $0 \leq x \leq z$. Due to Lemma \ref{lem:orth.pol} (vi) it holds that $\sup_{0 \leq x \leq z} x^\nu \phi_{\aopt-1}^2(x) \leq \nu^\nu |(\pol{\aopt-1}{1})'(0)|^{-\nu}$, $\nu \geq 0$.
The proof of Lemma 3.7 in \citet{Hanke} shows that \[
\polnorm{\pol{\aopt-1}{1}}{0}^{1/2}
	 \leq 
	\| \Proj_{z} \phi_{\aopt-1}(\KEst) \bEst\|_\calH.\] This yields with \ref{eq:source.H}
\begin{align*}
	&\|\KEst \fest{\aopt-1}{1} - \bEst\|_\calH
	 = 
	\polnorm{\pol{\aopt-1}{1}}{0}^{1/2}
	 \leq 
	\| \Proj_{z} \phi_{\aopt-1}(\KEst) \bEst\|_\calH\\
	& \leq 
	\| \Proj_{z} \phi_{\aopt-1}(\KEst) \KCov f\|_\calH
	+ 
	\| \Proj_{z} \phi_{\aopt-1}(\KEst) (\bEst-\KCov f)\|_\calH\\
	& \leq 
	\| \Proj_{z} \phi_{\aopt-1}(\KEst) \KCov f\|_\calH
	+ C_\epsilon  \gamma_n \left(\sup\limits_{0\leq x\leq z} \phi_{\aopt-1}^2\right)^{1/2}
	\\
	& \leq 
	\|\Proj_z \phi_{\aopt-1}(\KEst)\KEst^{\mu+1} u \|_\calH 
	+
	\|\Proj_z \phi_{\aopt-1}(\KEst)(\KEst^{\mu+1}-\KCov^{\mu+1}) u \|_\calH  + C_\epsilon\gamma_n\\
	& \leq 
	R\left\{
	 \left(\sup\limits_{0\leq x\leq z} x^{2\mu+2} \phi_{\aopt-1}^2\right)^{1/2}
	 +
	 C_\delta (\mu+1) \kappa^\mu \gamma_n 
	 \left(\sup\limits_{0\leq x\leq z} \phi_{\aopt-1}^2\right)^{1/2}
	\right\} + C_\epsilon\gamma_n\\
	& \leq 
	\pdiff{\aopt-1}^{-\mu-1} (2\mu+2)^{\mu+1} R + \{C_\delta (\mu+1) \kappa^\mu R  + C_\epsilon\}\gamma_n.
\end{align*}
This gives together with $\Ctau \gamma_n \leq \| \KEst \fest{\aopt-1}{1} - \bEst\|_\calH$
\begin{align*}
	\Ctau \gamma_n \leq \pdiff{\aopt-1}^{-\mu-1} (2\mu+2)^{\mu+1} R + \left\{C_\delta (\mu+1) \kappa^\mu R + C_\epsilon\right\}\gamma_n.
\end{align*}
If $\Ctau > 	C_\delta (\mu+1) \kappa^\mu R + C_\epsilon$ we finally have
\begin{equation}
\label{eq:diffm1bound}
	  \pdiff{\aopt-1} \leq \gamma_n^{-1/(\mu+1)} 
	  \left\{
	  	\frac{(2\mu+2)^{\mu+1} R}
	{
		\Ctau - 	C_\delta (\mu+1) \kappa^\mu R + C_\epsilon
	}
	  \right\}^{1/(\mu+1)}.
\end{equation}
If $\aopt=1$ it holds $\pol{\aopt-1}{1} = 1$ and thus $\pdiff{\aopt-1} = 0$ and the inequality (\ref{eq:diffm1bound}) is true as well.

We will derive an upper bound on $\left|
	\left(
		\pol{\aopt}{1}
	\right)'(0)
	-
	\left(
		\pol{\aopt-1}{1}
	\right)'(0)
\right|.$ Due to Corollary 2.6 of \citet{Hanke} we have
\begin{equation}
\label{eq:deriv.diff}
	\left|
	\left(
		\pol{\aopt-1}{1}
	\right)'(0)
	-
	\left(
		\pol{\aopt}{1}
	\right)'(0)
\right|\leq 
\frac{\polnorm{\pol{\aopt-1}{1}}{0}
}
{\polnorm{\pol{\aopt-1}{2}}{1}}.
\end{equation}
We have $0 \leq x \leq \xr{1}{\aopt-1}{1} < \xr{1}{\aopt-1}{2}$ due to the interlacing property of the roots in Lemma \ref{lem:orth.pol} (i) and thus $0 \leq \pol{\aopt-1}{2}(x) \leq 1$ for $0 \leq x \leq \xr{1}{\aopt-1}{2}$. 
With that we get with \ref{eq:source.H}
\begin{align*}
&\|\pol{\aopt-1}{1}(\KEst) \bEst \|_\calH 
\leq \polnorm{\pol{\aopt-1}{2}}{0}^{1/2} \\
&\leq 
\|\Proj_x \pol{\aopt-1}{2}(\KEst) \bEst\|_\calH + 
x^{-1/2}\|\Proj_x^\perp \KEst^{1/2} \pol{\aopt-1}{2}(\KEst) \bEst\|_\calH\\
&\leq 
\|\Proj_x \pol{\aopt-1}{2}(\KEst)(\bEst - \KCov f)\|_\calH + 
\|\Proj_x \pol{\aopt-1}{2}(\KEst) \KCov^{\mu+1} u\|_\calH + 
x^{-1/2}\polnorm{\pol{\aopt-1}{2}}{1}^{1/2}\\
&\leq 
C_\epsilon\gamma_n + R\left\{
C_\delta(\mu+1)\kappa^\mu \gamma_n
+
x^{\mu+1} 
\right\}
+
x^{-1/2}\polnorm{\pol{\aopt-1}{2}}{1}^{1/2}.
\end{align*}
For the choice $x_\ast = (C_x \gamma)^{1/(\mu+1)}$ we get
\[
	\polnorm{\pol{\aopt-1}{1}}{0}^{1/2} 
	\leq 
	\gamma_n\left\{
		C_\epsilon + C_\delta(\mu+1)\kappa^\mu R
		+C_x 
	\right\}	
+
x_\ast^{-1/2}\polnorm{\pol{\aopt-1}{2}}{1}^{1/2}.
\]
It holds $\polnorm{\pol{\aopt-1}{1}}{0}^{1/2} = \| \KEst \fest{\aopt-1}{1} - \bEst\|_\calH \geq \Ctau \gamma_n$. This yields with $\Ctau > C_\epsilon + C_x R + C_\delta (\mu+1)\kappa^\mu R$
\[
	\polnorm{\pol{\aopt-1}{1}}{0} \leq 
	\gamma_n^{-1/(\mu+1)} C_x^{-1/(\mu+1)}	
	\left\{
		1  - \frac{C_\epsilon + C_x R + C_\delta (\mu+1)\kappa^\mu R}{\Ctau}
	\right\}^{-2}
	\polnorm{\pol{\aopt-1}{2}}{1}.
\]
Together with (\ref{eq:deriv.diff}) we have
\[
	\left|
	\left(
		\pol{\aopt-1}{1}
	\right)'(0)
	-
	\left(
		\pol{\aopt}{1}
	\right)'(0)
\right|\leq 
\gamma^{-1/(\mu+1)}_n C_x^{-1/(\mu+1)}	
	\left\{
		1  - \frac{C_\epsilon + C_x R + C_\delta (\mu+1)\kappa^\mu R}{\Ctau}
	\right\}^{-2}.
\]
Combining this with (\ref{eq:diffm1bound}) completes the proof. \eop

\subsubsection{Proof of Theorem \ref{th:kpls}}
The proof is an application of Lemmas \ref{lem:error.difference} - \ref{lem:derivative} to (\ref{eq:main.inequality}). First note that $r \geq 3/2$ implies $\mu \geq 1$ and thus this condition in Lemma \ref{lem:secondlemma} holds.

Let us choose $x_\ast = (C_x \gamma_n)^{1/(\mu+1)}$. Lemma \ref{lem:orth.pol} (v) shows that $\pdifftwo{i}{r} = \sum_{j=1}^i (\xr{j}{i}{r})^{-1}$ for $i =1,\dots,n$, $r \in \N_0$.
Thus it holds $\pdifftwo{i}{1}^{-1} \leq \xr{1}{i}{1}$. 

Equation (\ref{eq:diffm1bound}) thus shows that $C_x$ can be chosen small enough such that
\[
	x_\ast \leq \pdifftwo{\aopt-1}{1}^{-1} \leq \xr{1}{\aopt-1}{1}
\]
and $C_x <1$, which makes the first condition in Lemma \ref{lem:error.difference} and \ref{lem:derivative} hold true. The choice $C = C_\epsilon + (\mu+1)\kappa^\mu R (1 + C_\delta)$ gives the second condition. 

Now we need to check the remaining condition of Lemma \ref{lem:secondlemma}, namely that a $C_z$ can be chosen such that $(C_z \gamma_n)^{1/(\mu+1)} \leq \xr{1}{\aopt}{1}$ is true. Lemma \ref{lem:derivative} yields a $C_z>0$ such that $C_z \gamma_n^{1/(\mu+1)} \leq \pdiff{\aopt}^{-1} \leq \xr{1}{\aopt}{1}.$ 
Denote $z_\ast = (C_z \gamma_n)^{1/(\mu+1)}$ and Lemma \ref{lem:secondlemma} can be applied. 

To ease notation we will denote everything in the derived bounds that does not depend on $\gamma_n$ as a constant $c_j$, $j \in \N$.
Thus we get by combining Lemmas \ref{lem:secondlemma} and \ref{lem:derivative} that with probability at least $1-\nu$
\begin{align*}
\|f - \fest{\aopt}{1}\|_\calH^2 &\leq c_1 \gamma_n^{\mu/(\mu+1)} + c_2 \gamma_n + c_3 \gamma_n^{1-1/(\mu+1)} + c_4 \gamma_n \pdiff{\aopt}\\
& \leq c_1 \gamma_n^{\mu/(\mu+1)} + c_2 \gamma_n + c_3 \gamma_n^{\mu/(\mu+1)} + c_5 \gamma_n^{1-1/(\mu+1)} = O\{\gamma_n^{\mu/(\mu+1)}\}
\end{align*}
and
\begin{align*}
	&\|\KEst^{1/2}(f - \fest{\aopt}{1})\|_\calH^2\\
	&\leq 
	c_6 \gamma_n^{(\mu+1/2)/(\mu+1)} +c_7 \gamma_n^{1/(2\mu+2)}\gamma_n + c_8 \gamma_n^{-1/(2\mu+2)}\gamma_n + c_9\gamma_n^{1/(2\mu+1)}\gamma_n \pdiff{\aopt}\\
	& \leq 
	c_6 \gamma_n^{(\mu+1/2)/(\mu+1)} +c_7 \gamma_n^{(2\mu+3)/(2\mu+2)} + c_8 \gamma_n^{(2\mu+1)/(2\mu+2)} + c_10\gamma_n^{1+1/(2\mu+2)-1/(\mu+1)}\\
	& = O\{\gamma_n^{(2\mu+1)/(2\mu+2)}\}.
\end{align*}
Finally Lemma \ref{lem:error.difference} gives
\[
	\|f_\aopt - \fest{\aopt}{1}\|_\calH^2 = O\{\gamma_n^{\mu/(\mu+1)}\}, ~~~~
	\|\KEst^{1/2}(f_\aopt - \fest{\aopt}{1})\|_\calH =O\{\gamma_n^{(2\mu+1)/(2\mu+2)}\}.
\]
Combining the above with (\ref{eq:main.inequality}) yields
\begin{align*}
\|f - f_\aopt\|_\calH^2 &= O\{\gamma_n^{\mu/(\mu+1)}\},\\
		\|f^\ast - f_\aopt\|_2^2 &= O\{\gamma^{1/2} \gamma_n^{\mu/(\mu+1)}\} + O\{\gamma_n^{(2\mu+1)/(2\mu+2)}\} = O\{\gamma_n^{(2\mu+1)/(2\mu+2)}\},
\end{align*}
completing the proof with $\mu = r-1/2$. \eop

\subsection{Proof of Theorem \ref{th:kpls2}}
The overall design of this proof is similar to the one of Theorem \ref{th:kpls} and makes heavy use of results obtained in \citet{Blan10b}. 
\subsubsection{Preparation for the proof}
The stopping index \ref{eq:stopping2} can be reformulated with $\mu = r-1/2$ as
\begin{equation}
\label{eq:alt.stop2}
	\aopt = \min\{1 \leq a \leq n: \|\KEst \fest{a}{1}- \bEst\|_\calH \leq C \zeta_n\},
\end{equation}
with $\zeta_n = \max\{\sqrt{\lambda_n \ed}\gamma_n,\lambda_n^{\mu+1}\}$. 

We will derive the result in a similar way to Theorem \ref{th:kpls}. First it holds due to (\ref{eq:main.inequality})
\begin{align}
	&\|f_\aopt - f^\ast\|_2 =\|\Sg (f_\aopt - f)\|_\calH \leq \|\Sg (f_\aopt - \fest{\aopt}{1})\|_\calH + \|\Sg (\fest{\aopt}{1} - f)\|_\calH \notag \\
	& \leq 
	C_\psi \lambda^{1/2} \|f_\aopt - \fest{\aopt}{1}\|_\calH + C_\psi\|\KEst^{1/2}(f_\aopt - \fest{\aopt}{1})\|_\calH  
	+ \|\Sg(\fest{\aopt}{1}-f)\|_\calH. 
	\label{eq:main.inequality2}
\end{align} 
Now we prove the analogue versions of Lemma \ref{lem:error.difference} -- \ref{lem:derivative}:

\begin{lemma}
\label{lem:error.difference2}
Let $x = C_x \lambda_n$, $C_x>0$ such that $0< x < x_{1,\aopt-1}^{[2]}$. Choose $C>\tilde{c}_2$, with $\tilde{c}_1 = R \max\{1,C_\psi^2,\mu \kappa^{\mu-1} C_\delta\}$ and $\tilde{c}_2 = 2\max\{C_\psi C_\epsilon \sqrt{C_x+1}, \tilde{c}_1 C_x (C_x^\mu+1)\}$. Then it holds 
\begin{align*}
	\|f_\aopt - \fest{\aopt}{1}\|_\calH 
 	&\leq 
	\frac{C^3}{C_x(C-\tilde{c}_2)^2}\lambda^{-1}_n \zeta_n,\\
\|S_n^{1/2}(f_\aopt - \fest{\aopt}{1})\|_\calH &\leq
\frac{C^2}{C_x^{1/2}(C-\tilde{c}_2)} \lambda^{-1/2}_n \zeta_n.
\end{align*}
\end{lemma}
Proof: According to the proof of Lemma \ref{lem:error.difference2} we can focus on the case $0< \aopt < n$. Furthermore we have due to (\ref{eq:est.difference})
\begin{equation}
\label{eq:est.difference2}
 	\|f_\aopt - \fest{\aopt}{1}\|_\calH 
 	\leq 
	\frac{\polnorm{\pol{\aopt-1}{2}}{0}}{\polnorm{\pol{\aopt-1}{2}}{1}}\|\pol{\aopt}{1}(\KEst) \bEst\|_\calH.
\end{equation}

Using Lemma A.3 in \citet{Blan10b} (and the first line of its proof) we have for $0 < x < x_{1,\aopt-1}^{[2]}$ with $\tilde{c}_1 = R \max\{1,C_\psi^2,\mu \kappa^{\mu-1} C_\delta\}$
\begin{align}
\label{eq:norm.difference}
\polnorm{\pol{\aopt-1}{2}}{0}^{1/2} 
\leq 
C_\psi C_\epsilon	 \sqrt{x + \lambda_n}\sqrt{d_{\lambda_n}}\gamma_n 
+ 
\tilde{c}_1 x \{x^\mu + Z_\mu(\lambda_n)\}
+ 
x^{-1/2} \polnorm{\pol{\aopt-1}{2}}{1}^{1/2}.
\end{align}
Here we define $Z_\mu(\lambda) = \lambda^\mu \mathbb{I}(\mu \leq 1) + \gamma_n \mathbb{I}(\mu>1)$. Note that under the assumptions of the theorem it holds $Z_\mu(\lambda_n) \leq \lambda_n^\mu$.

 Choosing $x = C_x \lambda_n$ yields in  (\ref{eq:norm.difference})
 \begin{align*}
 &\polnorm{\pol{\aopt-1}{2}}{0}^{1/2}\\
 &\leq 
 C_\psi C_\epsilon \sqrt{C_x + 1} \sqrt{\lambda_n d_{\lambda_n}}\gamma_n
  + 
  \tilde{c}_1 C_x(C_x^\mu+1)\lambda_n^{\mu+1}
  + 
  C_x^{-1/2} \lambda_n^{-1/2} \polnorm{\pol{\aopt-1}{2}}{1}^{1/2}\\
	&\leq 
	\tilde{c}_2\max\{\sqrt{\lambda_n d_{\lambda_n}} \gamma_n,\lambda_n^{\mu+1}\} + C_x^{-1/2} \lambda_n^{-1/2} \polnorm{\pol{\aopt-1}{2}}{1}^{1/2}\\
	& = 
	\tilde{c}_2 \zeta_n + C_x^{-1/2} \lambda_n^{-1/2} \polnorm{\pol{\aopt-1}{2}}{1}^{1/2}
	,
 \end{align*}
 with $\tilde{c}_2 = 2\max\{C_\psi C_\epsilon \sqrt{C_x+1}, \tilde{c}_1 C_x (C_x^\mu+1)\}$.  Due to the stopping condition (\ref{eq:alt.stop2}) we know that 
\[
\polnorm{\pol{\aopt-1}{2}}{0}^{1/2} 
\geq 
\polnorm{\pol{\aopt-1}{1}}{0}^{1/2}
 = 
 \|\KEst \fest{a}{1}- \bEst\|_\calH 
 \geq 
 C \zeta_n.
 \] 
 This gives
 \begin{equation}
 	\label{eq:norm.upperbound}
 	\polnorm{\pol{\aopt-1}{2}}{0}^{1/2} \leq  \frac{C}{\sqrt{C_x}(C-\tilde{c}_2)} \lambda_n^{-1/2}  \polnorm{\pol{\aopt-1}{2}}{1}^{1/2}.
 \end{equation}
 Plugging this into (\ref{eq:est.difference2}) yields together with the definition of the stopping index $\aopt$
\[
	\|f_\aopt - \fest{\aopt}{1}\|_\calH
 	\leq 
	\frac{C^3}{C_x(C-\tilde{c}_2)^2}\lambda^{-1}_n \zeta_n.
\]
 In a similar way we derive for the second case
\begin{equation*}
 	\|S_n^{1/2}(f_\aopt - \fest{\aopt}{1})\|_\calH = \frac{\polnorm{\pol{\aopt}{1}}{0}}{\polnorm{\pol{\aopt-1}{2}}{1}^{1/2}}
 	\leq 
 	\frac{\polnorm{\pol{\aopt-1}{2}}{0}^{1/2}}{\polnorm{\pol{\aopt-1}{2}}{1}^{1/2}}\polnorm{\pol{\aopt}{1}}{0}.
\end{equation*}
An application of (\ref{eq:norm.upperbound}) yields
 \[
 	\|S_n^{1/2}(f_\aopt - \fest{\aopt}{1})\|_\calH \leq 
 	\frac{C^2}{\sqrt{C_x}(C-\tilde{c}_2)} \lambda_n^{-1/2}\zeta_n.
 \]
 \eop

\begin{lemma}
\label{lem:secondlemma2}
Denote $\tilde{c}_1 = R \max\{1,C_\psi^2,\mu \kappa^{\mu-1} C_\delta\}$. 
For any $i=1,\dots,n$ and $0 < x < x_{1,i}^{[1]}$ we have under the conditions of the theorem
	\begin{align*}
	\|S^{1/2}(f_i^{[1]}-f)\|_\calH &\leq 
	C_\psi\left[C_\delta+\sqrt{2}C_\epsilon + \lambda_n\left\{C_\delta \pdiff{i} + \sqrt{2}C_\epsilon x^{-1}\right\}\right]\sqrt{d_{\lambda_n}}\gamma_n\\
	&+ 
	\tilde{c}_1 (\sqrt{x}+\sqrt{\lambda_n})(x^\mu+\lambda_n^\mu)
	+
	(1+\sqrt{x^{-1}\lambda_n})x^{-1/2}\|\KEst \fest{i}{1}-\bEst\|_\calH.
	\end{align*}
\end{lemma}
Proof: Follow the proof of Lemma A.2 in \citet{Blan10b}. Note that $Z_\mu(\lambda_n) \leq \lambda_n^\mu$.\eop

\begin{lemma}
\label{lem:derivative2}
Let $C>\max\{C_\psi C_\epsilon,\tilde{c}_2\}$, where $\tilde{c}_2$ is given in Lemma \ref{lem:error.difference2}. Choose $x = C_x \lambda_n$ such that $0 < x \leq x_{1,\aopt-1}^{[1]}$. Then there exists a constant $c^\ast>0$ such that
	\begin{align}
		\pdiff{a^\ast} \leq c^\ast \lambda^{-1}_n.
	\end{align}
\end{lemma}
Proof: In analogue to Lemma \ref{lem:derivative} we will first derive an upper bound on $\pdiff{\aopt-1}$. Lemma A.1 in \citet{Blan10b} yields
	\begin{align*}
	\|\KEst \fest{\aopt-1}{1} - \bEst \|_\calH 
	& \leq 
	R (2\mu+2)^{\mu+1}\max\{1,C_\psi^{2\mu}\} \pdiff{\aopt-1}^{-\mu-1}\\
	 &+ 2R \mu\kappa^{\mu-1} \max\{1,C_\delta,C_\psi^{2\mu},C_\psi^{2\mu} C_\delta\} \pdiff{\aopt-1}^{-1} Z_\mu(\lambda_n)\\
	 &+C_\epsilon C_\psi \left\{
	 	\pdiff{\aopt-1}^{-1/2} + \sqrt{\lambda_n}
	 \right\} \sqrt{d_{\lambda_n}}\gamma_n.
	\end{align*}
	Denote $\tilde{c}_3 = R (2\mu+2)^{\mu+1}\max\{1,C_\psi^{2\mu}\}$ and $\tilde{c}_4 = 2R \mu\kappa^{\mu-1} \max\{1,C_\delta,C_\psi^{2\mu},C_\psi^{2\mu} C_\delta\}$. 
	
The definition of $\aopt$ gives $C \zeta_n \leq \|\KEst \fest{\aopt-1}{1} - \bEst \|_\calH$. 
Combining both inequalities, setting $x = C_x \lambda_n$ and keeping $\sqrt{\lambda_n d_{\lambda_n}} \gamma_n \leq \zeta_n$ in mind gives
\begin{align*}
	(C-C_\psi C_\epsilon C_\lambda) \zeta_n 
	&\leq
	\tilde{c}_3 \pdiff{\aopt-1}^{-\mu-1}
	+
	\tilde{c}_4 \pdiff{\aopt-1}^{-1} \lambda_n^\mu\\
	&+
	C_\epsilon C_\psi \pdiff{\aopt-1}^{-1/2} \sqrt{d_{\lambda_n}}\gamma_n\\
	& \leq 3 \max\left\{
		\tilde{c}_3 \pdiff{\aopt-1}^{-\mu-1},
		\tilde{c}_4 \pdiff{\aopt-1}^{-1}\lambda_n^\mu,\right.\\
		&\left.
		C_\epsilon C_\psi \pdiff{\aopt-1}^{-1/2} \sqrt{d_{\lambda_n}} \gamma_n
	\right\}.
\end{align*}
Now we assume that the maximum on the right hand side is attained in each of the three possible cases
\begin{align*}
	\pdiff{\aopt-1} 
	& \leq 
	\{3 (C- C_\epsilon C_\psi)^{-1} \tilde{c}_3
	\}^{1/(\mu+1)} \zeta_n^{-1/(\mu+1)},\\
	\pdiff{\aopt-1} 
	& \leq 
	3 (C-C_\epsilon C_\psi)^{-1}\tilde{c}_4 \zeta_n^{-1} \lambda_n^\mu,\\
	\pdiff{\aopt-1} 
	& \leq 
	9 (C-C_\epsilon C_\psi)^{-2} C_\epsilon^2 C_\psi^2 \zeta_n^{-2} d_{\lambda_n} \gamma_n^2.
\end{align*}
Take $\tilde{c}_5 = \max[\{3 (C- C_\epsilon C_\psi)^{-1} \tilde{c}_3
	\}^{1/(\mu+1)},
	3 (C-C_\epsilon C_\psi)^{-1}\tilde{c}_4,
	9 (C-C_\epsilon C_\psi)^{-2} C_\epsilon^2 C_\psi^2]$.

It is easy to see that $\zeta_n^{-1/(\mu+1)}, \zeta_n^{-1}\lambda_n^\mu$ and $\zeta_n^{-2} d_{\lambda _n}\gamma_n^2$ are all bound from above by $\lambda_n^{-1}$. Hence we get 

\begin{equation}
\label{eq:first.derivative.bound}
\pdiff{\aopt-1}  \leq \tilde{c}_5 \lambda^{-1}_n.
\end{equation}
For the final step in the proof we have due to (\ref{eq:deriv.diff})
\begin{equation*}
	\left|
	\left(
		\pol{\aopt-1}{1}
	\right)'(0)
	-
	\left(
		\pol{\aopt}{1}
	\right)'(0)
\right|\leq 
\frac{\polnorm{\pol{\aopt-1}{1}}{0}
}
{\polnorm{\pol{\aopt-1}{2}}{1}}.
\end{equation*}
It holds $\polnorm{\pol{\aopt-1}{1}}{0} \leq \polnorm{\pol{\aopt-1}{2}}{0}$ and hence (\ref{eq:norm.upperbound}) yields
\begin{equation}
	\label{eq:first.derivative.bound2}
	\left|
	\left(
		\pol{\aopt-1}{1}
	\right)'(0)
	-
	\left(
		\pol{\aopt}{1}
	\right)'(0)
\right|\leq 
	\frac{C^2}{C_x(C-\tilde{c}_2)^2} \lambda^{-1}_n.
\end{equation}
The proof is complete by combining (\ref{eq:first.derivative.bound}) and (\ref{eq:first.derivative.bound2}).\eop 

\subsubsection{Proof of Theorem \ref{th:kpls2}}
We first restrict ourselves to the set where all concentration inequalities stated in the theorem hold simultaneously with probability at least $1-\nu$, $\nu \in (0,1]$. We only proof the convergence rates in the $\calL^2$-norm, the corresponding rates in the $\calH$-norm are done in the same way.
 
The theorem is proven by an application of Lemmas \ref{lem:error.difference2}--\ref{lem:derivative2}. To that end we need to check the conditions of those. Equation (\ref{eq:first.derivative.bound}) and the proof of Theorem \ref{th:kpls} show that we can take $C_x = \min\{1/2,\tilde{c}_5\}$ to fulfill $0 < x \leq x_{1,\aopt-1}^{[1]}$. Furthermore we can take $C=4 R \max\{1, C_\psi^2(\nu),(r-1/2)\kappa^{r-3/2} C_\delta(\nu),2^{-1/2}R^{-1} C_\psi(\nu) C_\epsilon(\nu)\}$ and the conditions of Lemma \ref{lem:error.difference2} and \ref{lem:derivative2} hold. Note that $x_{1,\aopt-1}^{[1]} \leq x_{1,\aopt-1}^{[2]}$ due to the interlacing property of the roots, see Lemma \ref{lem:orth.pol} (i).

For Lemma \ref{lem:secondlemma2} we need to find a $0 <  z < x_{1,\aopt}^{[1]}$. By Lemma \ref{lem:derivative2} there exists a constant $c^\ast>0$ such that
\[
	(c_\ast)^{-1} \lambda_n \leq \pdiff{\aopt}^{-1} \leq x_{1,\aopt}^{[1]},
\]
hence we choose $C_z = \min\{1/2,1/c^\ast\}$. Now, applying Lemmas \ref{lem:error.difference2}--\ref{lem:derivative2} to (\ref{eq:main.inequality2}) gives the result (we again denote any constant that does not depend on $n$ with $C_i$, $i \in \N$)
\begin{align*}
	&\|f_\aopt - f^\ast\|_2 \leq C_\psi \lambda_n^{1/2} \|f_\aopt - \fest{\aopt}{1}\|_\calH + C_\psi\|\KEst^{1/2}(f_\aopt - \fest{\aopt}{1})\|_\calH  
	+ \|\Sg(\fest{\aopt}{1}-f)\|_\calH \\
	 \leq &
	C_1 \lambda_n^{-1/2} \zeta_n   
	+ C_2 \lambda_n \pdiff{i} \sqrt{d_{\lambda_n}}\gamma_n
	+ C_3 \lambda_n^{\mu+1/2}
	+ C_4 \lambda_n^{-1/2} \|\KEst \fest{\aopt}{1}-\bEst\|_\calH.\\
	 \leq& 
	 C_5 \lambda_n^{-1/2} \zeta_n 
	 + C_6 \sqrt{d_{\lambda_n}}\gamma_n
	 + C_4 \lambda_n^{\mu+1/2} \leq \max\{C_4,C_5,C_6\} \lambda_n^{-1/2}\zeta_n.
\end{align*}
The error bound in the $\calH$-norm is proven in an analogue fashion.
\eop

\subsection{Proof of Corollary \ref{cor:pol.ed}}
Take $\lambda_n = \gamma_n^{2/(2r+s)}$. It is immediate that $\lambda_n^{r-1/2} = \gamma_n^{(2r-1)/(2r+s)} \geq \gamma_n$ for $n$ sufficiently large, hence the inequality (\ref{eq:lambda.inequ}) holds as soon as $\gamma_n \leq 1$. Then we have by Theorem \ref{th:kpls2} that 
\[
	\|f_{\alphaest_{a^\ast}} - f^\ast\|_2 = O\left\{\lambda_n^{-1/2}\zeta_n(\lambda_n)\right\} = O\left\{\gamma_n^{2r/(2r+s)}\right\}. 
\]
\eop

\subsection{Proof of Corollary \ref{cor:log.ed}}
Set $\lambda_n = \gamma_n^{1/r} \log\{1/(2r) \gamma_n^{-2}\}$. It is immediate that $\lambda_n \rightarrow 0$ as $\gamma_n$ converges to zero. For $r=1/2$ condition (\ref{eq:lambda.inequ}) holds trivially. Let $r>1/2$, then we have
\[
	\lambda_n^{r-1/2} = \gamma_n^{(r-1)/r} \log\{(r-1/2)/(2r) \gamma_n^{-2}\} \geq \gamma_n,
\]
This is equivalent to $2r-1 \geq 2 \exp(\gamma_n) \gamma_n^2$, which holds for $n$ sufficiently large and $r > 1/2$.

For the convergence rate we first show that $d_{\lambda_n} \gamma_n^2 \leq \lambda_n^{2r}$. We have
\[
	d_{\lambda_n} \gamma_n^2 = 
	\log\left\{
			1+ \frac{a}{\gamma_n^{1/r} \log^{1/r}(1/2 \gamma_n^{-1})}
		\right\}
	 \gamma_n^2 
	 \leq 
	 \log\left(
	\gamma_n^{-2}	
	\right)\gamma_n^2.
\]
Equivalently we need 
$a^r \leq \gamma_n (\gamma_n^{-2}-1)^r \log(1/2\gamma_n^{-2})$. As $\gamma_n$ converges to zero $\gamma_n(\gamma_n^{-2} -1)^{r}$ goes to infinity for any $r> 1/2$. Hence for suitably large $n$ it holds $\lambda_n^r \geq \sqrt{d_{\lambda_n}} \gamma_n$. Then the convergence rate is $\lambda_n^{-1/2} \zeta_n(\lambda_n) = \lambda_n^r = \gamma_n \log(1/2 \gamma^{-2})$. 

Because the convergence rate does not depend on $r \geq 1/2$ we can set $r = 1/2$.  \eop

\subsection{Proof of Theorem \ref{th:conc.equality}}
\subsubsection{Preparation for the proof}
We denote with $\mathrm{tr} (A)$ the trace of a trace class operator $A: \calH \rightarrow \calH$ and the tensor product $(f_1 \otimes f_2) h = \langle f_1, h \rangle_\calH f_2$ for functions $f_1,f_2,h \in \calH$. We use the notation $k_t = k(\cdot, X_t)$. Note that it holds $\|A\|_{\HS}^2 = \mathrm{tr}(A^\ast A)$ for a Hilbert-Schmidt operator $A$.

\begin{lemma}
\label{lem:sup1}
Under the assumptions \ref{con:k1} and \ref{con:k2} the following hold
\begin{enumerate}
\item[(i)] $\trace\{(k_t \otimes k_t) (k_s \otimes k_s)\} = k^2(X_t,X_s)$,
\item[(ii)] $\|\KCov\|_{\HS}^2 = \int_{\R^d}\int_{\R^d} k^2(x,y) \diff \Prob^{X_0}(x) \diff \Prob^{X_0}(y)$,
\item[(iii)] $\E [ \trace\{(k_0 \otimes k_0) \KCov\}] = \|\KCov\|_{\HS}^2$.
\item[(iv)] Let $X'$ and $X''$ be independent and identically distributed and denote $k' = k(\cdot,X')$, $k''=k(\cdot,X'')$. It holds for $\nu = 1,2$ and $\lambda>0$
\begin{align*}
	\mathrm{E}\left[
		\mathrm{tr}^\nu \left\{(S+\lambda)^{-1} k' \otimes k' k'' \otimes k''\right\}
	\right]
	= 
	\mathrm{tr}^\nu\{(S+\lambda)^{-1} S^2\}.
\end{align*}
\end{enumerate}
\end{lemma}
{\it Proof:} (i) Let $\{v_i\}_{i \in \N}$ denote an  orthonormal base of $\calH$. Then it holds due to the reproducing property (\ref{eq:rep.property})
\begin{align*}
	\mathrm{tr}\left\{(k_t \otimes k_t)(k_s \otimes k_s)\right\} = 
	\sum_{i=1}^\infty \langle v_i,k_t\rangle_\calH \langle v_i,k_s\rangle_\calH k(X_t,X_s) = \left\langle
		\sum_{i=1}^\infty \langle v_i, k_s \rangle_\calH v_i, k_t
	\right\rangle_\calH k(X_t,X_s).
\end{align*}

(ii)
\begin{align*}
\|\KCov\|^2_{\HS} &= \sum_{i=1}^\infty \langle \KCov v_i, \KCov v_i\rangle_\calH = \sum_{i=1}^\infty  \int\limits_{\R^d} \langle \KCov v_i, k(\cdot,x) \rangle_\calH \langle v_i, k(\cdot,x)\rangle_\calH \diff\Prob^{X}(x)\\
&= \int\limits_{\R^d} \int\limits_{\R^d} \left\langle \sum_{i=1}^\infty \langle v_i, k(\cdot,x) \rangle_\calH v_i, k(\cdot,y)\right\rangle_\calH k(x,y) \diff\Prob^{X}(x) \diff\Prob^{X}(y).
\end{align*}
The assertion follows because $\Prob^{X} = \Prob^{X_0}$.

(iii) 
\begin{align*}
		\E [\trace\{(k_0 \otimes k_0) \KCov \}] &= \E (\langle S k_0, k_0\rangle_\calH) = \E \left(\int\limits_{\R^d} \langle k_0, k(\cdot,x)\rangle^2_\calH \diff \Prob^{X}(x)\right) \\
		&= \int\limits_{\R^{d}}\int\limits_{\R^{d}} k^2(x,y) \diff\Prob^{X}(x) \diff\Prob^{X_0}(y) = \|S\|_{\HS}^2.
\end{align*}

(iv) Because $S$ is a compact operator the spectral decomposition $S = \sum_{i=1}^\infty \mu_i \psi_i \otimes \psi_i$ holds (recall that $\{\mu_i$, $\psi_i\}_{i=1}^\infty$ is the eigensystem of $S$). Let $k(\cdot,x) = \sum_{i=1}^\infty \alpha_i(x) \psi_i$, $x \in \R^d$. 
For $\nu=1$ we have
\begin{align*}
	\mathrm{E}\left[\mathrm{tr}\{(S+\lambda)^{-1} k' \otimes k' k''\otimes k''\}\right]& = 
	\E \left\{k(X',X'') \sum_{i=1}^\infty \langle \psi_i, k''\rangle_\calH \langle(S+\lambda)^{-1} k', \psi_i \rangle_\calH \right\}\\
	& =\sum\limits_{i=1}^\infty \frac{1}{\mu_i+\lambda} 
	 \E \left\{k(X',X'') \langle \psi_i, k''\rangle_\calH \langle  \psi_i,k' \rangle_\calH\right\}\\
	& = \sum\limits_{i=1}^\infty \frac{1}{\mu_i+\lambda} \sum\limits_{j=1}^\infty \E\left\{\alpha_j(X')\alpha_i(X')\alpha_j(X'')\alpha_i(X'')\right\}
\end{align*}
On the other hand
\begin{align*}
	\mathrm{tr}\left\{
		(S+\lambda)^{-1} S^2
	\right\} & = \mathrm{tr}\left\{
		(S+\lambda)^{-1} \E(k' \otimes k' k'' \otimes k'')
	\right\}\\
	& = \sum\limits_{i=1}^\infty \langle
		(S+\lambda)^{-1} \E \{k(X',X'') \langle \psi_i,k''\rangle_\calH k'\} 
	,\psi_i \rangle_\calH\\
	& = \sum\limits_{i=1}^\infty \frac{1}{\mu_i + \lambda }
	\langle
		\E \{k(X',X'') \langle \psi_i,k''\rangle_\calH k'\} 
	,\psi_i \rangle_\calH\\
	& = \sum\limits_{i=1}^\infty \frac{1}{\mu_i + \lambda }\sum\limits_{j=1}^\infty \E\{\alpha_j(X')\alpha_i(X')\alpha_j(X'') \alpha_i(X'')\}
\end{align*}
and we are done. The proof for $\nu=2$ is along the same lines. \eop

\begin{lemma}
\label{lem:acf.inequality}
Assume that condition \ref{D2} holds. Then we have
\begin{align}
\label{eq:sum}
	n^{-2} \sum\limits_{h=1}^{n-1}(n-h) |\rho_h| \leq C(q) \left\{
		\begin{array}{clc}
			n^{-1}&,& q>1\\
			n^{-1} \log(n)&,& q=1\\
			n^{-q}
			&,& q \in (0,1).
		\end{array}
	\right.,
\end{align}
with $C(q) = \zeta(q) \mathbb{I}(q>1) + \{5-\log(4)\}\mathbb{I}(q=1) + \{2(1-q)^{-1}-(2-q)^{-1}+(2-q)^{-1} 2^{2-q}\}\mathbb{I}\{q\in(0,1)\}$. Here $\zeta$ denotes the Riemann zeta function. 
\end{lemma}
{\textit  {Proof}}: Recall that by condition \ref{D2} we have $|\rho_h| \leq (h+1)^{-q}$, $h=0,\dots,n-1$ for some $q>0$. 

First assume $q \in (0,1]$. The integral test for series convergence gives lower and upper bounds for the hyperharmonic series as
\[
(1-q)^{-1}\{
	(n+1)^{1-q}-2^{1-q}
\} 
\leq 
\sum\limits_{h=2}^{n} h^{-q} 
\leq 2^{-q} 
+ 
(1-q)^{-1} 
\{
	n^{1-q}-2^{1-q}
\}.
\]
This yields
\begin{align}
&n^{-2} \sum\limits_{h=1}^{n-1} (n-h) (h+1)^{-q}  = 
n^{-2} \sum\limits_{h=2}^{n} (n+1-h)h^{-q} 
=
n^{-2} \left\{
	(n+1)\sum\limits_{h=2}^n h^{-q} - \sum\limits_{h=2}^n h^{-(q-1)}
\right\}
\notag \\
&\leq 
n^{-2} \left[
	(n+1) \left\{
		2^{-q} + (1-q)^{-1}(n^{1-q}-2^{1-q}
	\right\}
	-
	(2-q)^{-1}
	\left\{
		(n+1)^{2-q} - 2^{2-q}
	\right\}
\right].\label{eq:sum.inequality}
\end{align}
Now let $q\in (0,1)$, then it holds from (\ref{eq:sum.inequality}) and the fact that $n^{-2} \leq n^{-1} \leq n^{-q}$ 
\begin{align*}
	&n^{-2} \sum\limits_{h=1}^{n-1} (n-h) (h+1)^{-q}\\
	\leq & \frac{n+1}{n^2}\left\{
		\frac{2^{-q}(1-q)-2^{1-q}}{1-q}
	\right\}
	+ \frac{n+1}{n^{1+q}}(1-q)^{-1} - \frac{(n+1)^{2-q}}{n^2}(2-q)^{-1}
	+\frac{1}{n^2}\frac{2^{2-q}}{2-q}\\
	\leq & n^{-q}  [	
	\{
		2(1-q)^{-1} - (2-q)^{-1}
	\}
	+
	(2-q)^{-1}2^{2-q}
	],
\end{align*}
due to $2^{-q} (1-q) - 2^{1-q}<0$.

For $q=1$ we evaluate the limit 
\begin{align*}
	&\lim\limits_{q \rightarrow 1\pm} n^{-2} \left[
	(n+1) \left\{
		2^{-q} + (1-q)^{-1}(n^{1-q}-2^{1-q}
	\right\}
	-
	(2-q)^{-1}
	\left\{
		(n+1)^{2-q} - 2^{2-q}
	\right\}
\right]\\
 & = (2 n^2)^{-1} 
 [
 	3-\log(4) - n \{1+\log(4)\}
 ] + n^{-2} (n+1) \log(n)\\
 & \leq 
 \frac{\log(n)}{n} \left[
 	5 - \log(4)
 \right].
\end{align*}
The case $q>1$ is clear because the zeta-function $\zeta(q)$ is defined as the hyperharmonic series with coefficient $q$.\eop

Denote with $g_{h}$ the common density of $(X_h,X_0)^\T$ and $g_0$ the density of $X_0$. The next lemma and the subsequent corollary will be used to show that the quantities appearing in the sums of Theorem \ref{th:conc.equality} (i) can be linked to the autocorrelation function $\rho$:
\begin{lemma}
\label{lem:normal.quantity}
Under the assumptions \ref{con:k1}, \ref{con:k2} and \ref{D1} it holds for $h>0$ with $\rho_h = \tau_0^{-1} \tau_h$
\begin{align*}
	\int\limits_{\R^{2d}} k^2(x,y) \{g_h(x,y)-g_0(x)g_0(y)\}\diff(x,y) 
	&\leq 
	\frac{\kappa^2}
{\{(4 \pi \tau_0)^{d} \det(\NCov)\}^{1/2}}
\theta^{1/2}(\rho_h),\\
	\int\limits_{\R^{2d}} k(x,y)f(x)f(y) 
	\{g_h(x,y)-g_0(x)g_0(y)\}
	\diff(x,y) 
&\leq 
\frac{\kappa M}
{\{(4 \pi \tau_0)^{d} \det(\NCov)\}^{1/2}}
\theta^{1/2}(\rho_h),
\end{align*}
with $\theta(\rho) = 1 + (1- \rho^2)^{-d/2}
	- 	
	2^{d+1}
	(4 - \rho^2)^{-d/2}$, $\rho \in [0,1)$.
\end{lemma}
{\textit  {Proof}}: We will only proof the first inequality, the second one follows in the same way. 

By Jensen's inequality and \ref{con:k2} we know
\begin{align*}
	&\int\limits_{\R^{2d}} k^2(x,y) \{g_h(x,y)-g_0(x)g_0(y)\} \diff(x,y) \\
&\leq \kappa^2 \left[
	\int\limits_{\R^{2d}} \left\{
		g_h^2(x,y) - 2 g_h(x,y) g_0(x) g_0(y) + g_0^2(x) g_0^2(y) 
	\right\} \diff (x,y)
\right]^{1/2}.
\end{align*}
The first and third integral term can readily be calculated as
\begin{align*}
	\int\limits_{\R^{2d}} 
		g_h^2(x,y)\diff (x,y) 
		& = 
		[(4 \pi)^d (\tau^2_0-\tau^2_h)^{d/2} \det(\Sigma)]^{-1}\\
		\left\{\int\limits_{\R^{d}} 
		g_0^2(x)\diff x\right\}^2 
		& =
		\{(4 \pi)^d \tau^{d}_0 \det(\Sigma)\}^{-1}. 
\end{align*}
For the first equality we use $\det(A \otimes \Sigma) = \det(A)^d \det(\Sigma)^2$ for $A \in \R^{2 \times 2}$ and thus
\begin{equation}
\label{eq:normalsquared}
\int\limits_{\R^{2d}} 
		g_h(x,y) g_0(x) g_0(y)\diff (x,y) = 
		\frac{\int_{\R^{2d}}
			\exp\left(
				-1/2 z^\T G^{-1} z
			\right)		
		\diff z} 
		{(2\pi)^{2d} \det(\Sigma)^2\tau^d_0 (\tau^2_0-\tau^2_h)^{d/2}},
\end{equation}
with 
\[
	G^{-1} = 
	\left\{\left(
		\begin{matrix}
			\tau_0 & \tau_h\\
			\tau_h & \tau_0
		\end{matrix} 
	\right)^{-1}
	+
	\left(
			\begin{matrix}
			\tau^{-1}_0 & 0\\
			0 & \tau^{-1}_0
		\end{matrix} 
		\right)
	\right\}\otimes \Sigma^{-1}.
\]
It holds $\det(G) = 	(4 \tau^2_0 - \tau^2_h)^{-d} (\tau^4_0 - \tau^2_0 \tau^2_h)^d \det(\Sigma)^2$.
Thus we get with $(\ref{eq:normalsquared})$
\begin{align*}
		\int\limits_{\R^{2d}} 
		g_h(x,y) g_0(x) g_0(y)\diff (x,y) &= 
		\frac{(2\pi)^d 
			\tau^d_0 (\tau^2_0 - \tau^2_h)^{d/2} \det(\Sigma)		
		}
		{(2\pi)^{2d} \det(\Sigma)^2 (4 \tau^2_0 - \tau^2_h)^{d/2}\tau^d_0 (\tau^2_0-\tau^2_h)^{d/2}}\\
		& = 
		\left\{
			(2 \pi)^d 
			(
				4 \tau^2_0 - \tau^2_h
			)^{d/2}
			\det(\Sigma)
		\right\}^{-1},		
\end{align*}
completing the proof by multiplying all terms with $\tau^{-d}_0 \tau^d_0$.\eop

\begin{corollary}
\label{cor:rho.bound}
Under the assumptions \ref{con:k1}, \ref{con:k2}, \ref{D1} and \ref{D2} it holds for all $h>0$ and $q >0$
\begin{align*}
	\int\limits_{\R^{2d}} k^2(x,y) \{g_h(x,y)-g_0(x)g_0(y)\}\diff(x,y)  
	&\leq 
	\frac{
\kappa^2 d^{1/2}
}{
\{(2 \pi)^d \det(\Sigma)\}^{1/2}
} 
		(1-4^{-q})^{-1/4(d-2)}
|\rho_h|\\
	\int\limits_{\R^{2d}} k(x,y)f(x)f(y) \{g_h(x,y)-g_0(x)g_0(y)\}\diff(x,y) 
	&\leq 
	\frac{
\kappa M d^{1/2}
}{
\{(2 \pi)^d \det(\Sigma)\}^{1/2}
} 
		(1-4^{-q})^{-1/4(d-2)}
|\rho_h|.
\end{align*}
\end{corollary}
{\textit  {Proof}}: Recall that $\theta(\rho) = 1 + \{1-\rho^2\}^{-d/2} - 2^{d+1}\{4-\rho^2\}^{-d/2}$ for $\rho \in [0,1)$. We seek to find bounds on $\theta$ and the corollary can be proven by an application of Lemma \ref{lem:normal.quantity}. 

By assumption \ref{D2} we know there is a $\rho_\ast$ such that $\rho_h^2 \leq \rho_\ast^2<1$ for all $h > 0$. Thus consider $\rho \in [0,\rho_\ast]$. We start by finding a constant $C>0$ with
\[
	\theta'(\rho) = \,\rho \left\{
		(1-\rho^2)^{-d/2-1}
		-
		2^{d+1}
		(4-\rho^2)^{-d/2-1}
	\right\} d \leq  C \rho^2.
\]
Thus $C$ can be taken as $C = d \left\{
		(1-\rho_\ast^2)^{-d/2-1}
		-
		2^{d+1}
		(4-\rho_\ast^2)^{-d/2-1}
	\right\}$.

Thus we know that the slope of $\theta$ is always less than that of $C \rho^2$. Finally it holds that $\theta(0) = 0$ and thus $0 \leq \theta(\rho)\leq C \rho^2$, $\rho \in [0,\rho_\ast]$.

Under condition \ref{D2} it holds $\{1-\rho^2_\ast\}^{-d/2} \leq \{1-2^{-2 q}\}^{-d/2}$, completing the proof by using Lemma \ref{lem:normal.quantity}.
\eop	

The final preparatory result is used to derive the probablistic bound in Theorem \ref{th:conc.equality} (iv) and is similar to Corollary \ref{cor:rho.bound}:
\begin{lemma}
	\label{lem:normal.quantity2}
	Under the assumptions \ref{con:k1}, \ref{con:k2}, \ref{D1} and \ref{D2} it holds for $\lambda>0$
\begin{align*}
	 \int\limits_{\R^{2d}} k(x,y) \langle (S+\lambda)^{-1} k(\cdot,x),k(\cdot,y)\rangle_\calH
  \left\{ 
 f_h(x,y)-f_0(x)f_0(y)
 \right\}\diff (x,y)
 \leq 
	\tilde{c} |\rho_h|\ed 	,
\end{align*}	
with $\tilde{c} =  \sqrt{d \{1-4^{-q}\}^{-d-1}} \kappa$.
\end{lemma}
Proof: Denote $\beta(x,y) = \langle (S+\lambda)^{-1} k(\cdot,x),k(\cdot,y)\rangle_\calH$. By the Cauchy-Schwarz inequality
\begin{align}
\notag
	\phi_h 
	&= 
	\int\limits_{\R^{2d}} k(x,y) \beta(x,y)
 \left\{ 
 f_h(x,y)-f_0(x)f_0(y)
 \right\}\diff (x,y)\\
 \label{eq:two.integral.cauchy}
 & \leq  
\left[\int\limits_{\R^{2d}} k^2(x,y) \beta^2(x,y) f_0(x) f_0(y)\diff (x,y) \int\limits_{\R^{2d}}\left\{
\frac{f_h(x,y)}{\sqrt{f_0(x)f_0(y)}}-\sqrt{f_0(x)f_0(y)}
\right\}^2\diff (x,y)\right]^{1/2}.
\end{align}
Denote by $X'$ and $X''$ two independent copies of $X_0$ and $k' = k(\cdot,X')$, $k''=k(\cdot,X'')$.
We start by bounding the first integral term in the product:
\begin{align*}
\int\limits_{\R^{2d}} k^2(x,y) \beta^2(x,y) f_0(x) f_0(y)\diff (x,y)
& = 
\E\left\{
	k^2(X',X'')\langle (S+\lambda)^{-1} k',k''\rangle_\calH^2\right\}\\
	& = \E\left[
		\mathrm{tr}^2\left\{(S+\lambda)^{-1} k' \otimes k' k'' \otimes k''
	\right\}
\right]\\
&\leq \kappa^2 \mathrm{tr}^2\{(S+\lambda)^{-1} S\} = \kappa^2 \ed^2.
\end{align*}
In the second to last inequality we used Lemma \ref{lem:sup1} (iv) and the definition of $\ed$.

The second integral in the product in (\ref{eq:two.integral.cauchy}) is
\begin{align*}
	\int\limits_{\R^{2d}}\left\{
\frac{f_h(x,y)}{\sqrt{f_0(x)f_0(y)}}-\sqrt{f_0(x)f_0(y)}
\right\}^2\diff (x,y) & =
\int\limits_{\R^{2d}}
\frac{f_h(x,y)}{\sqrt{f_0(x)f_0(y)}}\diff (x,y) -1.
\end{align*}
We proceed in the same way as in the proof of Lemma \ref{lem:normal.quantity} by using properties of the Gaussian distributions at hand.

First we have 
\[
	F_h(x,y) = \frac{f_h^2(x,y)}{f_0(x)f_0(y)} = \frac{(2\pi)^d \tau_0^d \mathrm{det}(\Sigma)}{(2\pi)^{2d} \mathrm{det(\Sigma_h)}} \exp\left[
		-\frac{1}{2} (x^\T,y^\T)
		\left\{
			2\Sigma_h^{-1} - \frac{1}{\tau_0} \left(
			\begin{matrix}
				\Sigma & 0\\
				0 & \Sigma
			\end{matrix}
			\right) 
		\right\}^{-1} \left(
		\begin{matrix}
		x\\
		y
		\end{matrix}
		\right)
	\right].
\]
Denote $G^{-1} = 2\Sigma_h^{-1} - \tau_0^{-1} \left(
			\begin{matrix}
				\Sigma & 0\\
				0 & \Sigma
			\end{matrix}
			\right) $. It follows in a similar way as in the proof of Lemma \ref{lem:normal.quantity} that $\mathrm{det}(G) = \tau_0^{2d} \mathrm{det}^2(\Sigma)$ and $\mathrm{det}(\Sigma_h) = (\tau_0^2 - \tau_h^2)^d \mathrm{det}^2(\Sigma)$. Hence we have
\begin{align*}
	\int\limits_{\R^{2d}} F_h(x,y) \diff (x,y) &= 
	\frac{(2\pi)^d \tau_0^d \mathrm{det}(\Sigma)}{(2\pi)^{2d}(\tau_0^2 - \tau_h^2)^d \mathrm{det}^2(\Sigma)} (2\pi)^d \tau_0^d \mathrm{det}(\Sigma)\\
	& = 
	\frac{\tau_0^{2d}}{(\tau_0^2-\tau_h^2)^d} 
	= 
	\frac{1}{(1-\rho_h^2)^d}.
\end{align*}
Under \ref{D2} we have $\rho_h < 1$ for all $h>0$ and there is a $\tilde{\rho} = \max_h |\rho_h| \leq 2^{-q} < 1$ and hence it holds $(1-\rho_h^2)^{-d} \leq d (1-4^{-q})^{-d-1} \rho_h^2$ and we are done. \eop

\subsubsection{Proof of the theorem}
First note that the the operator norm is dominated by the Hilbert-Schmidt norm.
By Markov's inequality we have for $\nu \in (0,1]$
\begin{align*}
	\Prob\left(
		\|\KEst - \KCov\|_\HS^2 \leq \nu^{-1} \E\|\KEst - \KCov\|_\HS^2
	\right) &\geq 1 - \nu,\\
	\Prob\left(
		\|\bEst - \KCov f\|_\calH^2 \leq \nu^{-1} \E\|\bEst - \KCov f\|_\calH^2
	\right) &\geq 1 - \nu.
\end{align*}
(i) It holds due to $\KEst = n^{-1} \sum_{t=1}^n k_t \otimes k_t$
\begin{align*}
	\E\left(\|\KEst - \KCov\|^2_{\HS}\right) = \frac{1}{n^{2}}\sum\limits_{t,s=1}^n \left(
		\E [\trace\{(k_t \otimes k_t)(k_s \otimes k_s)\}] - 2 \E [\trace\{( k_0 \otimes k_0) \KCov\}] + \|\KCov\|^2_{\HS}
	\right).
\end{align*}
For the first summand we get $\E [\trace\{( k_t \otimes k_t)( k_s \otimes k_s)\}] = \E \{k^2(X_t,X_s)\}$, due to Lemma \ref{lem:sup1} (i). 
Using the stationarity of $\{X_t\}_{t=1}^n$ and Lemma \ref{lem:sup1} (iii) we get
\begin{align*}
	\E\left(\|\KEst- \KCov\|^2_{\HS}\right) =  \frac{1}{n} \left\{\E \{k^2(X_0,X_0)\} - \|\KCov\|^2_{\HS} \right\} 
	+  \frac{2}{n^{2}} \sum\limits_{h=1}^{n-1} (n-h) \left[
		\E \{k^2(X_h,X_0)\} - \|\KCov\|_{\HS}^2
	\right],
\end{align*}
yielding the first result by an application of Lemma \ref{lem:sup1} (ii).\ \\
For the second equation we see due to the independence of $\{X_t\}_{t=1}^n$ and $\{\varepsilon_t\}_{t=1}^n$ that
\[
	\|T^\ast_n y - \KCov f\|^2_\calH = \sigma^2 n^{-1} \E \{k(X_0,X_0)\} + \E \left(\|\KEst f - \KCov f\|^2_\calH\right).
\]
The rest follows along the same lines as the first part of the proof.

(ii) 
An application of part (i) of this theorem, Corollary \ref{cor:rho.bound} and Lemma \ref{lem:acf.inequality} yields this result.

(iii) Because the $\{\varepsilon_t\}_{t\in \Z}$ are independent and identically distributed and $\{X_t\}_{t \in \Z}$ is stationary it holds
\begin{align*}
	&\mathrm{E}\left\{\left\|
	(S+\lambda)^{-1/2}(S_n f - T_n^\ast y)
	\right\|^2_\calH\right\}\\
	& = 
	\mathrm{E}\left\{\left\|
		(S+\lambda)^{-1/2} T_n^\ast \varepsilon
	\right\|_\calH^2\right\} 
	=
	n^{-2} \sum\limits_{t,s=1}^n \mathrm{E}\left\{
	\left\langle
		\varepsilon_t (S+\lambda)^{-1} k_t, \varepsilon_s k_s
	\right\rangle_\calH\right\}\\
	& = 
	n^{-1}\sigma^2  \mathrm{E}\left\{
	\left\langle
		(S+\lambda)^{-1} k_0,  k_0
	\right\rangle_\calH \right\}
	= 
	n^{-1}\sigma^2  \mathrm{E}\left\{
	\left\|(S+\lambda)^{-1/2} k_0\right\|^2_\calH\right\}.
\end{align*}
By the definition of $\ed$ we get 
\[
	\mathrm{E}\{\|(S+\lambda)^{-1/2}k_0\|_\calH^2\} 
= 
\mathrm{E}[\mathrm{tr}\{
	(S+\lambda)^{-1} k_0 \otimes k_0
\}]
= \mathrm{tr}\{
	(S+\lambda)^{-1} S
\} 
= 
\ed.
\]
Using $n^{-1/2} \leq \gamma_n(q)$ proves the result.

(iv) Consider first
\begin{align*}
\E\left\{ \|(S+\lambda)^{-1/2}(\KEst - S)\|_\HS^2\right\} &=n^{-2}\sum_{t,s=1}^n \E\left[ \mathrm{tr}\{ (S+\lambda)^{-1}(k_t \otimes k_t - S)(k_s \otimes k_s - S)\}\right]\\
&=n^{-1} \E\|(S+\lambda)^{-1/2}(k_0 \otimes k_0 - S)\|_\HS^2\\
& + 
\sum_{h=1}^{n-1} \E\left[
	\mathrm{tr}\{(S+\lambda)^{-1}(k_0 \otimes k_0 - S)(k_h \otimes k_h - S)\}
\right].
\end{align*}
Continuing with the expression inside the sums we expand
\begin{align*}
	\phi_h &= \E\left[
	\mathrm{tr}\{(S+\lambda)^{-1}(k_0 \otimes k_0 - S)(k_h \otimes k_h - S)\}
\right]\\
	&= \E\left[
	\mathrm{tr}\{(S+\lambda)^{-1}
	k_0 \otimes k_0k_h \otimes k_h\}\right] -  \mathrm{tr}\{(S+\lambda)^{-1}S^2\}\\
	&=
	 \E\left\{
		k(X_0,X_h) \langle (S+\lambda)^{-1} k_0, k_h\rangle_\calH
	\right\}
	-
	\mathrm{tr}\{(S+\lambda)^{-1}S^2\}.
\end{align*}
Using Lemma \ref{lem:sup1} (iv) we see that 
\[
	\mathrm{tr}\{(S+\lambda)^{-1} S^2\} = \int\limits_{\R^{2d}} k(x,y) \langle (S+\lambda)^{-1} k(\cdot,x),k(\cdot,y)\rangle_\calH \mathrm{d}\Prob^X(x) \mathrm{d}\Prob^X(y).
\]
Hence we have
\[
 \phi_h = \int\limits_{\R^{2d}} k(x,y) \langle (S+\lambda)^{-1} k(\cdot,x),k(\cdot,y)\rangle_\calH
 \left\{ 
 \mathrm{d}\Prob^{X_0,X_h}(x,y)-
 \mathrm{d}\Prob^{X_0}(x) \mathrm{d}\Prob^{X_0}(y)\right\}.
\]
This can be bound by the results of Lemma \ref{lem:normal.quantity2} and together with Lemma \ref{lem:acf.inequality} there exists a constant $C(q)>0$ such that with probability at least $1-\nu$
\[
	\|(S+\lambda)^{-1/2}(\KEst - S)\|_\calL^2 \leq \nu^{-1} C^2(q) \gamma^2_n(q) \ed,
\]
with $\gamma^2_n(q) = \left\{\begin{array}{cc}
	n^{-1},& q>1\\
	n^{-1} \log(n), & q=1\\
	n^{-q}, & q\in (0,1).
\end{array}
\right.$

This implies $\|(S+\lambda)^{-1/2}(\KEst - S)(S+\lambda)^{-1/2}\|_\calL \leq \nu^{-1/2} C(q) \lambda^{-1/2}\sqrt{\ed}\gamma_n(q)$. Let $\lambda = \lambda_n$ be a sequence converging to zero such that $\lambda_n^{-1/2}\sqrt{d_{\lambda_n}} \gamma_n(q) \rightarrow 0$. Let $n$ be large enough such that $\nu^{-1/2} C(q)\lambda_n^{-1/2} \sqrt{d_{\lambda_n}} \gamma_n(q) < 1$. Using Lemma A.5 in \citet{Blan10b} we obtain 
\[
	\|(S+\lambda)^{1/2}(\KEst + \lambda)^{-1/2}\| \leq  [1-\nu^{-1/2} C(q) \lambda_n^{-1/2} \sqrt{d_{\lambda_n}} \gamma_n(q)]^{-1/2} \leq \sqrt{2}.
	\] 
	The latter inequality can be fulfilled for $n$ large enough such that $\nu^{-1/2} C(q) \lambda_n^{-1/2} \sqrt{d_{\lambda_n}} \gamma_n(q) \leq 1/2$.\eop

\subsection{Proof of Proposition \ref{prop:source}}

Recall that $\KCov u = \E \{u(X_0) k(\cdot,X_0)\}$ for $u \in \calH$. Define the independent random variables $Y_1,\dots,Y_\mu$ that are all distributed as $X_0$.

First consider the following observation for $\mu \in \N$:
\begin{align}
	S^\mu u &= S (S^{\mu-1}u) = \E_{Y_1} \{(S^{\mu-1}u)(Y_1) k(\cdot,Y_1)\} =\E_{Y_2}\E_{Y_1}\{(S^{\mu-2}u)(Y_2) k(Y_1,Y_2) k(\cdot,Y_1)\} = \dots \notag \\
	 & = 
	\E_{Y_\mu} \cdots \E_{Y_1} \left\{\prod\limits_{\nu=1}^{\mu-1} k(Y_\nu,Y_{\nu+1}) u(Y_\mu) k(\cdot,Y_1)\right\}.
	\label{eq:iterated.expectation}
\end{align}
We take $u = \sum_{i=1}^\infty c_i k(\cdot,z_i)$ for $\{z_i\}_{i \in \N}, \{c_i\}_{i \in \N} \subset \R$ such that $\|u\|_\calH^2 = \sum_{i,j=1}^\infty c_i c_j k(z_i,z_j) \leq R^2$. The fact that a function $u \in \calH$ can be represented as a linear combination of kernel functions is due to the Moore-Aronszajn Theorem, see \citet{Berlinet}.

Define the matrix $\Gamma = [\Gamma_{i,j}]_{i,j=1}^{\mu+2} \in \R^{(\mu+2) \times (\mu+2)}$ via
\begin{equation*}
	\Gamma_{i,j} = \left\{
		\begin{array}{clc}
			\sigma_x^{-2} + 2 l &,& i=j = 2,\dots,\mu+1\\
			l &, & i=j=1,\mu+2\\
			-l &,& |i-j|=1\\
			0 &,& else
		\end{array}.
	\right.
\end{equation*}
Then we have via the integration of Gaussian functions and (\ref{eq:iterated.expectation})
\begin{align*}
	f(x) &= \frac{1}{(2\pi \sigma^2_x)^{\mu/2}} \sum\limits_{i=1}^\infty c_i \int\limits_{\R^\mu} \exp\left\{-1/2
		(x,x_1,\dots,x_\mu,z_i) \Gamma (x,x_1,\dots,x_\mu,z_i)^\T 
	\right\} \diff (x_1,\dots,x_\mu)\\
	& = 
	\frac{1}{\sigma^\mu_x\det(\Gamma_{2:\mu+1})^{1/2}} \sum\limits_{i=1}^\infty c_i 
	\exp\left[ - 1/2 \det(\Gamma_{2:\mu+1})^{-1}\left\{
		\det(\Lambda_{1:\mu+1}) (x^2 + z_i^2) - 2 l^{\mu+1} x z_i		
	\right\}\right].
\end{align*}
Here we used the symmetry property $\det(\Gamma_{2:\mu+2}) = \det(\Gamma_{1:\mu+1})$ as the first and last rows and columns of $\Gamma$ are identical.
This concludes the proof. \eop

\subsection{Proof of Proposition \ref{prop:source}}
In \citet{Shi08} it was shown that the eigenvalues of $S$ have the form $\mu_i = a b^{i-1}$, $i=1,2,\dots$ with
\begin{equation*}
	a = \sqrt{2}(1+\beta+\sqrt{1+\beta})^{-1/2},~~
	b = (1+\beta+\sqrt{1+2\beta})^{-1} \beta.
\end{equation*}
and $\beta = 4 l \sigma^2_x$. It is clear that $0<b<1$ and hence $0 < \mu_i \leq a$. We have $d_\lambda = \sum_{i=0}^\infty \{1+ a^{-1} b^{-i}\lambda\}^{-1}$. Denote $f(x) =\{1+ a^{-1} b^{-x}\lambda\}^{-1}$. We want to apply the integral test to the sum. We have $\int_0^\infty f(x) \mathrm{d}x = \log^{-1}(b^{-1})\log(1+a \lambda^{-1})$. This yields the bounds
\begin{align*}
	\frac{\log(1+a/\lambda)}{\log(b^{-1})} \leq 
	d_\lambda
	\leq 
	\frac{1}{1+\lambda/a} + \frac{\log(1+a/\lambda)}{\log(b^{-1})}.
\end{align*}
On $\lambda \in (0,1]$ we get $d_\lambda \leq D \log(1+a/\lambda)$ for a constant $D>0$. This can be seen as follows: The function $g_1(\lambda)=(1+\lambda/a)^{-1}$ is bounded from above by $C_1 = 1$ and the function $g_2(\lambda) = (b^{-1})\log(1+a/\lambda)$ is lower bounded by $c_2 = \log(1+a)$ and has no upper bound. 

Hence on the set $I=\{\lambda \in (0,1]: g_2(\lambda) \geq C_1\}$ we can choose $C = 2$. On the set $I^c$ we have on the other hand $C g_2(x) \geq c_2 C \geq g_1(x) + g_2(x)$, hence we need $C = 2 c_2^{-1} C_1 = 2 \log^{-1}(1+a)$. The choice $D = 2 \log^{-1}(b^{-1}) \max\{1, \log^{-1}(1+a)\}$ is sufficient and we have $d_\lambda \leq D  \log(1+a/\lambda)$, $\lambda \in (0,1]$.
\eop
\end{document}